\documentclass[numbers=noenddot,article,12pt]{scrreprt}

\usepackage[ngerman]{babel}	
\usepackage[latin1]{inputenc}
\usepackage{array}
\usepackage{amsmath}
\usepackage{amssymb}
\usepackage{latexsym}
\usepackage{textcomp}
\usepackage{wasysym}
\usepackage{theorem}
\usepackage{epsfig}
\usepackage{psfrag}
\usepackage{graphicx}
\usepackage{scrtime}
\usepackage{longtable}

\usepackage[all]{xy}
\usepackage{yfonts}
\usepackage{nicefrac}
\usepackage{stmaryrd}
\usepackage{multirow}
\usepackage{bigdelim}
\usepackage{tabularx}
\usepackage{amsfonts}

\makeatletter
\renewcommand\section{\@startsection 
    {section}{1}{0mm}
    {0.7\baselineskip}
    {0.35\baselineskip}
    {\bfseries\large\boldmath}
    } 
\makeatother 
\makeatletter
    \renewcommand\subsection{\@startsection
       {subsection}{2}{0mm}
       {0\baselineskip}
       {0\baselineskip}
       {\bfseries\normalsize\boldmath}
       }
\makeatother

\addtokomafont{chapter}{\rmfamily\bfseries\boldmath\Large}
\addtokomafont{chapter}{\Large}

  {\noindent\textbf{Beweis #1 : }}%
  {
  {\unskip\nobreak\hfill\penalty 50\hskip .001pt\mbox{}\nobreak\hfill\mbox{$\square$}\parfillskip=0pt\finalhyphendemerits=0\medbreak}\rm%
  }
  {\noindent\textbf{Beweisidee #1 : }}%
  {
  {\unskip\nobreak\hfill\penalty 50\hskip .001pt\mbox{}\nobreak\hfill\mbox{$\square$}\parfillskip=0pt\finalhyphendemerits=0\medbreak}\rm%
  }
\makeatother

\newtheorem{definition}{Definition}[section]
\newtheorem{proposition}[definition]{Proposition}
\newtheorem{lemma}[definition]{Lemma}

\newtheorem{remark}[definition]{Remark}
\newtheorem{beispiel*}{Beispiel}

\newtheorem{theorem}[definition]{Theorem}
\newtheorem{corollary}[definition]{Corollary}

\renewcommand{\epsilon}{\varepsilon}

\newcolumntype{P}[1]{>{\raggedright\arraybackslash}p{#1}}

\newdimen\headrulewidth
\newdimen\footrulewidth

\begin{document}

\begin{center}
\vspace{5cm}

\textbf{\Large Groups of piecewise isometric permutations}\\ 
\textbf{\Large of lattice points}\\ 
\textit{\Large by}\\ 
\textit{\Large Robert Bieri and Heike Sach}
\end{center}


\begin{small}
\textbf{Summary.} Let $M$ denote either Euclidean or hyperbolic $N$-space, and let $\Gamma$ be a discrete group of isometries of $M$, with the
property that  $\Gamma$ respects and acts tile-transitively on a convex polyhedral tesselation of $M$. Given an arbitrary base point $p \in M$, we
consider the orbit  $\Omega:= \Gamma p \subseteq M$ and define a notion of ``$\Gamma$-polyhedral pieces'' $S \subseteq \Omega$. The objects of our
interest are the groups $G_{\Gamma}(S)$ of \textit{all piecewise} $\Gamma$-\textit{isometric permutations} on S.
\\In this paper we merely present the two most basic examples, and these play rather different roles: The case when
$\Gamma = {\text{PSL}}_2(\mathbb Z)$ acting on the hyperbolic plane reveals that the groups $G_{\Gamma}(\Omega)$ here have
prominent relatives: they are closely related to Richard Thompson's group $V$. And in the case  $\Gamma =$ Isom$(\mathbb Z^N)$ we find that the groups $G_{\Gamma}(S)$ have diverse but to some extent
computable finiteness properties.
\\The conjunction of the two examples suggests that to investigate the piecewise $\Gamma$-isometric permutation groups more systematically might be a
worthwhile project and might yield interesting new groups with an accessible finiteness pattern.

\begin{footnotesize}
finiteness properties of groups, piecewise isometric infinite permutations.  
\end{footnotesize}
\end{small} 

\vspace{10mm}

\chapter{Introduction}
\section{Generalities and main result}
\subsection {The groups.}~
Let $M$ denote either Euclidean or hyperbolic $N$-space, $N\in \mathbb N$, and \mbox{$\Gamma \le$ Isom$(M)$} a discrete group of isometries of $M$
with the property that $\Gamma$ admits a finite sided convex fundamental polyhedron $D$ with finite volume\footnote{In the hyperbolic case this implies
that $D$ is actually a generalized polytope - see Thm. 6.4.8 in [Ra94].}. We aim to study certain groups of permutations of the orbit $\Omega: =
\Gamma p$, for a given point $p \in M$.

To define the notion of a piecewise $\Gamma$-isometric permutation $\pi: \Omega \to \Omega$ requires a notion  of ``$\Gamma$-polyhedral pieces'' of 
$\Omega$ on which $\pi$ should be isometric, and it is reasonable to require that the geometry of these pieces be related
to the geometry of $\Gamma$.
Thus, together with the base point $p \in M$ we choose a finite set $\mathcal{H}$ of ``$\Gamma$-\textit{relevant}'' \textit{closed
half-spaces} of $M$, and the resulting groups will - to some extent - depend on this choice. One possibility would be to write 
$D$ as the intersection $D = \bigcap_{H\in \mathcal{H}}H$ with $\mathcal{H}$ an irredundant finite set of half spaces, whose boundaries are spanned by the sides of $D$; but other choices might be more convenient.
\\
By a \textit{convex $\Gamma$-polyhedral} subset $P$ of $M$ we mean any finite intersection of 
$\Gamma$-translates H$\gamma$, where $\gamma \in\Gamma$ and $H\in \mathcal{H}$. And a general $\Gamma$-\textit{polyhedral} subset of $M$ is a finite 
union of convex ones. By abuse of language, we call the intersection $S = \Omega\cap P$ a \textit{(convex)} $\Gamma$-\textit{polyhedral piece of} 
$\Omega$ whenever $P\subseteq M$ is a (convex) $\Gamma$-polyhedral subset. However, we will also meet situations where it is reasonable to use
``$\Gamma$-convex pieces'' to define piecewise $\Gamma^\ast$-isometric  permutations on $\Omega$ when \mbox{$\Gamma^\ast
\le$ Isom$(M)$} contains $\Gamma$ as a subgroup and $\Gamma^\ast p = \Gamma p$.

\textbf{Definition.}  Let $S\subseteq\Omega$ be a $\Gamma$-polyhedral piece. A permutation $g: S \to S$ is said to be \emph{piecewise
$\Gamma^\ast$-isometric} if $S$ can be written as a union of finitely many convex $\Gamma$-polyhedral pieces $\Omega = S_1 \cup S_2 \cup... \cup S_k$ 
with the property that the restriction of $g$ to each $S_i$  is also the restriction of an isometry $\varphi \in \Gamma^\ast.$

We write $G_{\Gamma^\ast}(S)$ for the group of all piecewise isometric $\Gamma^\ast$-permutations of $S$. The permutations in $G_{\Gamma^\ast}(S)$ with finite support form a normal subgroup of $G_{\Gamma^\ast}(S)$ which we denote by $\Sigma_{\infty}(S)$; the quotient group $G_{\Gamma^\ast}(S)/\Sigma_{\infty}(S)$ is often particularly interesting.

In this paper we consider the group $G_{\Gamma^\ast}(S)$ in two special cases:

a) When $M$ is the hyperbolic plane and $\Gamma = {\text{PSL}}_2(\mathbb Z)$ we
consider the group phi$(\Omega):=  G_\Gamma(\Omega)$ on the orbit $\Omega = \Gamma e^{i\pi/3}$. We show that phi$(\Omega)$ is the group of all quasi
isomorphisms of the planar dyadic tree, and the quotient phi$(\Omega)/\Sigma_{\infty}(S)$ is Richard Thompson's group $V$. Thus, the
quotient $G_\Gamma(S)/\Sigma_\infty(S$) for $S$ an arbitrary hyperbolic or Euclidean polyhedral piece can be viewed as a far reaching generalization
of Thompson's group $V$.

b) Our main concern then is the case when $M$ is Euclidean $N$-space, $\Gamma$ is the translation group $T = $~Tra$(\mathbb Z^N)$, and $\Gamma^\ast$
is either the full isometry group Isom$(\mathbb Z^N)$ or equal to $T$.  In both cases we use $T$-polyhedral pieces $S\subseteq\mathbb Z^N$ and term
them the \emph{orthohedral} subsets of $\mathbb Z^N$. We consider the \textit{piecewise Euclidean isometry} groups pei$(S) = G_{\Gamma^\ast}(S)$ and 
the \textit{piecewise Euclidean translation} groups pet$(S) = G_{\Gamma}(S)$ of arbitrary orthohedral subsets $S\subseteq \mathbb Z^N.$ If $S$ is the 
set of all points on the positive coordinate axes then the pet-group $G_{\Gamma}(S)$ is the Houghton group $H_N$ \cite{ho78}. Hence the quotients
pei$(S)/\Sigma_{\infty}(S)$ can be viewed as Euclidean relatives of Thompson's group $V$.
\\We are interested in structural and finiteness properties of these generalized Houghton groups.

\subsection{The finiteness length of a group.}~
By definition, every group is \textit{of type} $F_0$; every finitely generated group is \textit{of type} $F_1$; every finitely presented group
(equivalently: every fundamental group $\pi_1(X)$ of a finite cell complex $X$) is \textit{of type} $F_2$ ; and $\pi_1(X)$ is \textit{of type} $F_m$
($m\geq2$) if $X$ is a finite cell complex and $\pi_1(X) = 0$, for all $i$ with $2\leq i<m$.

Ten years after C.T.C. Wall introduced these finiteness properties, Borel and Serre \cite{bs73}/\cite{bs76} showed that all semi-simple
S-arithmetic groups have special homological features; in particular they are of type $F_{\infty}$ (i.e., type $F_m$ for all $m \in \mathbb N$). And
this was only the first of a number of important infinite families of groups that turned out to be of type $F_{\infty}$ in the following decades; many of them, just like arithmetic
groups, in the center of mainstream group theory: automorphism groups of free groups \cite{cv86}, Richard Thompson' s groups \cite{bg84}, etc.
More recent results in this direction are based on Ken Brown's topological discrete Morse theory technique \cite{br87} and its powerful
CAT$(0)$-version of Bestvina-Brady \cite{bb97}.

The insight that many important groups have much further reaching finiteness properties than finite presentability is great progress - but having ``good'' finiteness properties is only one 
side of the concept: The focus on the  \textit{finiteness length} function $fl: \mathbf{Gr}
\rightarrow \mathbb N\cup \{0,\infty\}$, defined on all groups $G$ by
\[fl(G) := \{\sup{m~|~G~ \mbox{is of type}~ F_m}\}\]
takes both sides into account. Analogous \textit{algebraic length functions} $afl_A$ are defined for every $G$-module $A$,
to be the supremum of all non-negative integers $m$ with the property that $A$ admits a free resolution which is finitely
generated in  all dimensions $\leq m$. The functions $afl_A$  have the considerable advantage, that they extend
immediately to monoids $G$. We write $afl$ for $afl_\mathbb Z$, where $\mathbb Z$ stands for the infinite cyclic up with the trivial
$G$-action; by the Hurewicz Theorem we know that $afl$ coincides with $fl$ on all finitely presented groups (i.e.,
whenever $fl(G) \geq 2)$. An important feature of both $fl$ and $afl_A$ is that they are constant on 
commensurability classes of groups.

The finiteness length of a stray group is notoriously difficult to compute. Nevertheless, to study and interpret
accessible parts of the pattern that these functions carve into group theory can be very fruitful. A convincing example is
the following: If we fix a finitely generated group $G$, then the function $Hom(G,\mathbb R_{add}) \rightarrow \mathbb N 
\cup\{0,\infty\}$, which associates with each homomorphism $\chi: G \rightarrow \mathbb R_{add}$ the value of $afl_A$ on the 
submonoid $\chi^{-1}([0,\infty)) \subseteq G$, imposes in the finite dimensional $\mathbb R$-vector space $Hom(G,\mathbb
R_{add})$ the pattern exhibited by the homological $\Sigma$-invariants $\Sigma^k(G;A)$ of \cite{br88}. On the other hand, we can 
also evaluate $fl$ and $afl_A$ on the commensurability classes of subgroups containing $G'$, and this yields patterns on
the rational Grassman space of $\mathbb Q$-linear subspaces of $G/G'\otimes\mathbb Q$ (which parametrizes these classes).
The core of the main $\Sigma$-results of \cite{bns86}, \cite{br88}, \cite{r89}, \cite{bge03} consists then of exhibiting
the precise relationship between the two patterns. 
\\An intriguing point is that in all computable examples, the finiteness length patterns have a polyhedral flavor: they
turn out to be expressible in terms of finitely many inequalities. One of the few general results here, polyhedrality of 
$\Sigma^0(G;A)$ when $G$ is Abelian, was proved in \cite{bgr84} by methods which were later partly re-detected in tropical
geometry. But polyhedrality questions on $\Sigma^k(G;A$) for non-Abelian $G$ and $k>0$ are wide open.

\subsection{The results.}~ 
In this paper we make first steps to evaluate the finiteness length function $fl$ and $afl$ on what we like
to view as the \emph{pei- and pet-clouds around} Isom$(\mathbb Z^N)$, resp. $\mathbb Z^N$: the groups pei$(S)$, resp.
pet$(S)$, as $S$ runs through all orthohedral subsets of $\mathbb Z^N$.
Our main tool is Ken Brown's approach in \cite{br87}.
In that paper each $fl$-result comes together with a parallel $afl$-result, and hence our results have the same feature.

To state our main results requires the following notation: By an \emph{orthant of rank} $n$ ($n\in
\mathbb N$) we mean any subset $L\subseteq\mathbb Z^N$ isometric to the standard rank-$n$ orthant $\mathbb N^n$. Each
orthohedral set $S\subseteq\mathbb Z^N$ is the disjoint union of finitely many orthants $S = L_1\cup L_2\cup \ldots\cup L_k.$ By the \emph{rank} of $S$, denoted $\text{rk}S$, we mean the maximum rank of the orthants $L_i$; and the \emph{height} of $S$, denoted $h(S)$, is the number of orthants of rank $\text{rk}S$ among the $L_i$.

\textbf{Theorem A.} \emph{Rank and height of an orthohedral set S determines the group \emph{pei}$(S)$ up to isomorphism, and we have
\mbox{$fl\big( $\emph{pei}$(S)\big) \ge h(S) - 1.$}}

Fore a more precise result see Theorem 3.1. As $h(\mathbb Z^n) = 2^n$ we have, in particular,

\textbf{Corollary.} $fl\big(\text{pei}(\mathbb Z^n )\big) \ge2^n - 1.$

The isomorphism class of the pet-groups are not determined by rank and height of the orthohedral set $S$. Here we find a 
precise result for the special case when $S$ is a stack of $h$ (parallel) orthants of the same rank:

\textbf{Theorem B.} \emph{If $S$ is a stack of $h$ rank-n orthants then fl\big(\emph{pet}$(S)\big) = h(S) - 1.$}

For a generalization to a stack of $k$-skeletons of an orthant see Proposition 4.2 in Section 4.1.

\subsection{Outlook.}~
Let $\Gamma$ be a discrete group of (Euclidean or hyperbolic) isometries with polyhedral fundamental domain of finite volume. By generalizing the
definition of the group pei$(\mathbb Z^N)$ to the groups pi$_\Gamma(\Omega)$ of all piecewise $\Gamma$-isometric permutations of the orbit $\Omega =
\Gamma p$, we have endowed each such group $\Gamma$ with the pi$_\Gamma$-\emph{cloud} of all piecewise $\Gamma$-isometric permutation groups
pi$_\Gamma(S)$ where $S$ runs through the $\Gamma$-polyhedral subsets of $\Omega$.  Our success with evaluating the finiteness length function on the clouds around
Isom$(\mathbb Z^N)$ and $\mathbb Z^N$, together with the observation that the groups around $SL_2(\mathbb Z)$ are closely related to the highly 
respected Thompson groups, indicates that aiming to investigate the finiteness pattern on more of these
pi$_\Gamma$-\emph{clouds} might be a difficult but worthwhile program.

\subsection{Remark on the history of this paper.}~
C. Houghton introduced his groups in \cite{ho78}. Theorem B, in the Hougton
group case, i.e., when S is a stack of rays, is due to K.S. Brown \cite{br87}.
The inequality $fl$\big(pet$(S)\big) \ge h(S) - 1$ in the special case when $S$ is a stack of
quadrants (as well as the equality for a certain ``diagonal subgroup'' of pet$(S)$) is due to the second author and
appears in her Diploma thesis (Frankfurt 1992 \cite{sa92}), to which the first author contributed little more then the
definition of the groups. It was recent increasing interest for the Houghton groups (\cite{abm14}, \cite{bcmr14}, \cite{st15}, \cite{leh08},
\cite{za15} etc.) that suggested  Heike Sach's Diploma thesis \cite{sa92} should be published, translated, and generalized.

The surprising observation that our concept leads to rather interesting connections in the hyperbolic case is more recent.

\section{A hyperbolic example}
\subsection{$\mathbf{{\text{PSL}}_2(\mathbb Z)}$-polyhedral subsets of the hyperbolic
plane.}~
Let $M$ denote the upper half-plane model of the hyperbolic plane $\mathbf{H}^2$ on which the group  $\Gamma =
{\text{PSL}}_2(\mathbb Z)$ acts by Moebius transformations.
Its standard fundamental polyhedron is the intersection $D = H\cap H'\cap H'' \subseteq \mathbf{H}^2$ of the three half spaces 
\[H = \{z\in\mathbb C \mid |z|\ge 1\},  H' = \{z\in\mathbb C \mid \mbox{Re}(z) \le \frac{1}{2}\},  H'' = \{z\in\mathbb C \mid
- \frac{1}{2} \le \mbox{Re}(z)\}.\]

$D$ is a generalized triangle with one finite edge which we denote by $e$. The union $ T := \bigcup_{\gamma\in\Gamma}\gamma e$ is the \emph{Serre
tree}, a combinatorial planar dyadic tree on which $\Gamma$ acts by planar tree-isomorphisms. In fact, restriction to $T$ induces an isomorphism between $\Gamma$
and the group of all \emph{planar automorphisms} of $T$ i.e., the tree automorphism which respect the cyclic orientation of the link of each vertex.

The situation here is particularly simple: $\Gamma$ acts transitively on the oriented edges of $T$, hence the half-planes $H, H', H''$ are
$\Gamma$-translates of each other. Moreover, the union of all $\Gamma$-translates of $\partial H$ is the union of all $\Gamma$-translates of $\partial D$. Hence all $2$-dimensional convex
$\Gamma$-polyhedral subsets $P\subseteq\mathbf{H}^2$ are unions of $\Gamma$-translates of tiles $D\gamma , \gamma \in \Gamma.$
Since each $\Gamma$-translate of $\partial H$ intersects $T$ in an edge, we observe that the intersection $T\cap\partial P$ of the Serre tree $T$ with 
the boundary of an arbitrary convex $\Gamma$-polyhedral set $P \subseteq \mathbf{H}^2$ is a finite set of tree-geodesic segments $[v,v']$ of $T$, and
$T\cap P$ is a sub-forest of $T$ in which all but finitely many vertices are of degree 3.

In terms of the horo-ball $B = \{z\in\mathbb C~|~\mbox{Im}(z) \ge 1\}$ and its
$\Gamma$-translates we can be more precise: Each connected component $[v,v']$ of $T\cap\partial P$ consists of edges tangent to one of these horo-balls, and at its endpoints $v, v' \in \partial P$ turns into a ray that 
plunges vertically into either the same or a the neighboring horo-ball and runs to infinity. Each connected component of  $T\cap Int(P)$ is a 
homogeneous rooted tree hanging at an endpoint of a connected component of  $T\cap\partial P$.
 
\subsection{Piecewise $\mathbf{\Gamma}$-isometric versus piece planar tree-isometric permutations.}

Let $\Omega = \Gamma p$ be the set of all vertices of the Serre tree $T$. On the basis of the above description it easy to pin down the $\Gamma$-polyhedral
pieces of $\Omega$.
We need the following terminology concerning a subgraph  $X \subseteq T$.  A vertex $v$ of $X$ is \emph{inner} if its degree $\text{deg}_{X}(v) = 3$
and $v$ is a \emph{leaf} of $X$ if $\text{deg}_{X}(v) = 1. ~ X$ is a \emph{rooted subtree} of $T$ if all except one of its vertices are inner; the one
exceptional vertex $p$ is the \emph{root} of $X$. We write $T_0(p)$ for a rooted tree with root $p$ of degree $2$ and $T_1(p)$ if the root $p$ is a leaf.
Thus we have
\begin{proposition}\label{proposition2.1}
Each convex $\Gamma$-polyhedral piece $S \in\Omega$ has a cofinite subset, which is the disjoint union of the vertex sets
of finitely many rooted subtrees of the form $T_1(p)$.~~~~~~~~
\hfill$\square$
\end{proposition}

It follows that every piecewise $\Gamma$-isometric permutation of $\Omega =
\text{ver}(T)$ can be interpreted as a \emph{piecewise planar tree isometric}
(ppti) permutation of the vertices of the tree $T$ -- i.e, a permutation of ver$(T)$ that respect all but finitely many
edges and the link-orientations at all but finitely many vertices of $T$. Other authors use the term quasi-autmorphisms
(\cite{leh08}, \cite{ls07}, \cite{bmn13}, and \cite{nst15}).
Conversely, the convex closure of $T_1(p)$ is a convex $\Gamma$-polyhedral subset $P$ with $\text{ver}(T_1(p)) = \Omega\cap P$, and
the tree-isometric embedding $T_1(p)\rightarrow T$ extends to an isometric embedding of the convex closure
of $T_1(p)$ into $\mathbf{H}^2$. Hence we can summarize:

\begin{theorem}
If $\Gamma = {\emph{PSL}}(\mathbb Z)$ and $\Omega = \Gamma e^{i\pi / 3}$ is the set of vertices of the Serre tree $T$ then
\emph{phi}$(\Omega)$ coincides with \emph{ppti}$(\Omega)$. \hfill$\square$
\end{theorem}

\subsection{Classification of the $\mathbf{\Gamma}$-polyhedral subsets $\mathbf{S\subseteq \Omega}$.}~
To classify the $\Gamma$-polyhedral subsets $S\subseteq \Omega$ (i.e., the vertex set of finite unions of finite and
rooted subtrees of $T$) up to ppti-isomorphism is not too difficult. First of all, adjoining or removing finitely many edges and displacing connected 
components does not change the ppti-isomorphism type -- hence $T$ itself is
ppti-isomorphic to the disjoint union $T_0(p)\cup T_0(p')$.
\\
The observation that the disjoint union $T_0(p)\cup T_0(p')\cup {\{p''\}}$ is
ppti-isomorphic to $T_0(p'')$ leads to the following classification: Every $\Gamma$-polyhedral subsets $S\in \Omega$ is ppti-isomorphic to one of the following types: We define isomorphism types 
$T_k$ for each $k\in \mathbb{Z}$ as follows:
$T_0$ and $T_1$ stand for the types represented by single rooted sub-trees as above, and for $k\in \mathbb{N}$ we write 
$T_k$ for the type represented by the disjoint union of $T_0$ with $k$ isolated vertices, while $T_{-k}$  stands for
the forest obtained from $T_0$ by removing $k$ vertices. Note that the disjoint union of $k$ copies of $T_0$ is $T_{-k}$, 
in particular $T$ is phi-isomorphic to $T_{-1}$.

It should not be too difficult to show that the isomorphism types $T_k$ are pairwise different. The isomorphism types of
their quasi-automorphism groups $QT_k = \text{ppti}(\text{ver}(T_k))$ and their finiteness length seems to be more of a
problem! Brita Nucinkis and Simon St. John-Green \cite{nst15} have recently shown that $fl(QT_1) = \infty$, and (based on
\cite{leh08} and \cite{bmn13}) they uncovered structural properties which show, in particular, that the two
groups $QT_0$ and $QT_1$ cannot be isomorphic. It would be interesting to know more about the finiteness length and the
isomorphism types of the groups $QT_k$  for all integers $k$.
\vskip 2mm 
\subsection{Relationship with Thompson's group V.}~
We do not claim originality for the context of this subsection -- the quasi-automorphism group of dyadic trees and the 
$SL_2(\mathbb Z)$-aspects of Thompson's groups have been around for some time (e.g. \cite{ls07}, \cite{leh08}, 
\cite{fks11}).
\\
We start by choosing a finite subtree $X\subseteq T$ with the property that the given ppti-permutation $\alpha$ of
ver$(X)$ respects all edges and all link orientations outside $X$.
By adjoining edges we may and will assume that $X$ has no vertex of degree $2$. We write $\partial X$ for the set of leaves
of $X$ and $\mathring{X}$ for the subtree generated by the the inner vertices of $X$. The complement of $\mathring{X}$ in
ver$(T)$ contains for each leaf $a\in \partial X$ a leafless rooted tree $B_a$ with $B_a\cap X = \{a\};$ and we have $T  =  X\cup(\bigcup_{ a\in \partial X}B_a).$

The image of $\bigcup_{ a\in \partial X}B_a$ under $\alpha$ is the disjoint union of the leafless rooted subtrees $\alpha(B_a)$, $a\in\partial X$.
Let $X'$ denote the convex hull of $\alpha(\partial X)$. $X'$ is a finite subtree of $T$ which contains from each of the trees $B_a$ the root $a$,
and only the root $a$. And these roots $a$ are the leaves of $X'$. A vertex $b$ of degree $2$ of $X'$ would be the base point of a geodesic ray $R$ of
$T$ emanating from $b$ into the complement of $\alpha(\bigcup_{ a\in \partial X}B_a)$. On the one hand $R$ could enter any of the subtrees $\alpha(B_a)$
only through its root which is impossible; on the other hand $R$ cannot stay in the complement of $\alpha(\bigcup_{ a\in \partial X}B_a)$ since this
is finite.
Hence $X'$ contains no vertices of rank $2$, and this shows \mbox{$X' \cup (\alpha(\bigcup_{ a\in \partial X}B_a)) = T$}.
Thus ver$(X') = \alpha(\text{ver}(X))$,  $X'$ has the same features as $X$, and $\alpha |_{\text{ver}(X)}$ defines bijections of both $\partial
X \rightarrow \partial X'$ and $\text{ver}(\mathring{X}) \rightarrow \text{ver}(\mathring{X}')$.

We summarize:
\begin{proposition}\label{proposition1.4} The \emph{ppti}-permutations of \emph{ver}$(T)$ are classified by the pairs $[X,X']$ 
of finite subtrees of $T$ which have no vertex of degree $2$, together with bijections
\par
\hspace{35mm}$	\alpha: \partial X\to\partial X'    ~~~~and~~~~ 	\mathring{\alpha}: \emph{ver} (\mathring{X})\to
\emph{ver}(\mathring{X}').$\hfill$\square$
\end{proposition}
This is closely related to the often used combinatorial description of Richard Thompson's group $V$ (see \cite{br87}), 
which is commonly defined in terms of the rooted tree $T_0(p)$; it is the group of homeomorphisms induced by the
almost-planar-automorphisms  $\alpha$ of $T_0(p)$ with $\alpha(p) = p$ on the Cantor set at the boundary at infinity of 
$T_0(p)$. The observation that $T$ is ppti-isomorphic to  the disjoint union of two copies of $T_0(p)$ and hence to 
$T_0(p) - \{p\}$ shows that $V$ is isomorphic to the group induced by ppti(ver$(T)$) on the boundary at infinity of $T$.
The kernel of the corresponding homomorphism \mbox{ppti$(\Omega$) $\twoheadrightarrow$ $V$} is the
normal subgroup $\Sigma_\infty(\Omega)$ consisting of all permutations of $\Omega$ with finite support.
Combined with Theorem 2.2 this proves 

\begin{theorem}
\emph{phi}$(\Omega)/\Sigma_\infty(\Omega)$ is isomorphic to Thompson's group $V$.
\end{theorem}

\chapter{The orthogonal-Euclidean case}

\section{Orthohedral sets}

\subsection{Integral orthants in $\mathbb Z^N $.}~ In the standard $N$-dimensional
Euclidean integral lattice $\mathbb Z^N $, endowed with the canonical basis $X$,
we consider affine-orthogonal transformations \mbox{$\tau_{a,A}: \mathbb Z^N \to
\mathbb Z^N$}, $\tau_{a,A}(x) = a+Ax$, where $A\in O(N,\mathbb Z)$ is an
integral orthogonal matrix and $a \in \mathbb Z^N $.
Inside $\mathbb Z^N$ we have the \emph{standard orthant of rank} $N, \mathbb N^N
\subseteq\mathbb Z^N$, and all images of its $k$-dimensional faces, $0\le k\le
N$, under affine-orthogonal transformations.
More precisely: the subsets $L = \tau_{a,A} \langle Y \rangle \subseteq\mathbb Z^N $, where $\langle Y \rangle $ stands for the monoid generated by the $k$-element subset $Y\subseteq X$.
We call $L$ an \emph{integral orthant} (of rank $k$, and based at $a\in L$) of
$\mathbb Z^N $ or just a \emph{rank-$k$} orthant. 

We write $\Omega^k$ for the set of all rank-$k$ orthants of $Z^N$ and $\Omega^\ast$ for the union $\bigcup_k\Omega^k$.
$\Omega^\ast$ is partially ordered by inclusion, with $\Omega^0=Z^N$. The subset of all orthants based at the origin 0
will be denoted by $\Omega_0^\ast \subseteq \Omega^\ast$; it retracts the order preserving projection $\tau:\Omega^\ast
\rightarrow \Omega_0^\ast$ which associates to each orthant $L\in \Omega^\ast$ based at $a\in\mathbb Z^N $ its unique
parallel translate $\tau(L)= -a + L\in\Omega_0^\ast$. $\tau(L)$ is characterized by its canonical basis $Y = \{y \in
\pm X~|~ a+\mathbb Ny \subseteq L\}$ which indicates the directions of $L$; hence we call $\tau(L)$ the indicator of $L$.
Note that $Y$ is given by the function $f$: $X$ $ \rightarrow \{0, 1, -1\}$ with  $f(x) = \epsilon\in \{1, -1\}$ if  
$\epsilon x\in Y$, and $f(x) = 0$ if  $\{x, -x\}\cap Y = \emptyset$; hence $|\Omega_0^\ast| = 3^N$. 

We call a subset $S\subseteq\mathbb Z^N$ \emph{orthohedral} if it is the union of a finite set of orthants - without loosing
generality we can assume that the union is disjoint. The \emph{rank} of $S$, denoted $\text{rk}S$, is the maximum rank of
an orthant contained in $S$. If $S$ is isometric to $\mathbb N^k \times \{1, 2, \ldots , h\}$, we call it a \emph{stack of orthants of height} $h$. The terminology agrees with the \emph{height} $h(S)$
of an arbitrary orthohedral set $S\subseteq\mathbb Z^N$, defined as the number of orthants of maximal rank, $\text{rk}S$, which
participate in a pairwise disjoint finite decomposition of $S=L_1\cup L_2\cup \ldots \cup L_m$ -- see §1.3 in the
Introduction.

We will often use the elementary

\begin{lemma}\label{lemma2.1}
Orthohedrality of subsets $S\subseteq\mathbb Z^N$ is closed under the set theoretic operations of taking intersections, 
unions, and complements.\hfill$\square$
\end{lemma}

We write $\Omega^k(S)=\{L\in\Omega^k | L\subseteq S\}$ for the set of all rank-$k$ orthants of $S$, $\Omega^\ast(S)$ for the disjoint union over k, and $\Omega_0^\ast(S)$ for the
set of all orthants of S based at the origin 0. We consider the restriction of the indicator map $\tau:\Omega^\ast(S)
\rightarrow \Omega_0^\ast$. We write $S_\tau\subseteq \mathbb Z^N$ for the union of all orthants in $\tau(\Omega^\ast(S))$
and call this the indicator image of S. Note that $\tau(\Omega^\ast(S))=\Omega_0^\ast(S_\tau)$, and we can view the
indicator map as a rank preserving surjection \mbox{$\tau:\Omega^\ast(S) \twoheadrightarrow \Omega_0^\ast(S_\tau)$}. 

\subsection{Germs of orthants.}~ Two orthants $L, L'$  in $\Omega^\ast$ are said to be
\emph{commensurable} if $\text{rk}L = \text{rk}(L\cap L') = \text{rk}L'.$ We write $\gamma (L)$ for the commensurability class of $L$ and call 
it the \emph{germ} of $L.$ The union of all members of $\gamma (L)$ is the coset of $\mathbb Z^N$ and denoted by
$\langle L\rangle$ $\subseteq\mathbb Z^N$. The germs inherit from their representing orthants $L$ the \emph{rank}, 
relations like \emph{parallelism and orthogonality}, and also a \emph{partial ordering} defined as follows: given two
germs $\gamma, \gamma'$ we put $\gamma\leq\gamma'$ if they can be represented by orthants $L, L' \in\Omega^\ast$ with
$L\subseteq L'$. Note that if $L$ and $L'$ are arbitrary orthants representing $\gamma$ and $\gamma'$, respectively, then
$\gamma\leq\gamma'$ if and only if (1) $L'$ contains an orthant parallel to $L$ and (2) $L\subseteq\langle L\rangle$.

We write $\Gamma^\ast(S)= \bigcup_k\Gamma^k(S)$ for the set of all germs of orthants in $S$ and $\Gamma_0^\ast(S)$ for
the set of all germs represented by an orthant of $S$ based at the origin 0. $\Gamma^\ast(\mathbb Z^N)$ and
$\Gamma_0^\ast(\mathbb Z^N)$ are abbreviated by  $\Gamma^\ast$ and  $\Gamma_0^\ast$, respectively. As $\Gamma_0^\ast$
and $\Omega_0^\ast$ are canonically bijective, we will identify them when this is convenient. Note that $\Gamma^\ast(S)$
is a convex subset of $\Gamma^\ast$ in the sense that if $\gamma\in\Gamma^\ast(S)$ then $\{\gamma'\in\Gamma^\ast~|~
\gamma'\leq\gamma\}\subseteq\Gamma^\ast(S)$. We can interpret the indicator map as an order and rank preserving surjection
\mbox{$\tau:\Gamma^\ast(S) \rightarrow \Gamma_0^\ast$} with $\tau(\Gamma^\ast(S))= \Gamma_0^\ast(S_\tau)$. By
max$\Gamma^\ast(S)$ we mean the set of all maximal germs of $S$.

\textbf{Excercise:} Observe that $\tau(\text{max}\Gamma^\ast(S)) \supseteq \text{max}\Gamma_0^\ast(S_\tau)$, but this is
not, in general, an equality. 

\begin{lemma}\label{lemma2.2} \emph{max}$\Gamma^\ast(S)$ is finite for each orthohedral set $S$. The set of all germs of
rank $n = {\emph rk}S$ is a subset of \emph{max}$\Gamma^\ast(S)$, whose cardinality coincides with the height $h(S)$.
Hence $h(S)$ is independent of the particular decomposition of $S$.
\end{lemma}
\textit{Proof.} Let $S = \bigcup_j L_j$ be an arbitrary decomposition of $S$ as a finite pairwise disjoint union of
orthants $L_j$. Each orthant $L\subseteq$ is the disjoint union of the orthants $M_j = L\cap L_j$, and
exactly one of them is commensurable to $L$. Hence $\gamma(L)  =  \gamma(M_j)\subseteq\gamma(L_j)$. This shows that each
germ $\gamma\in\Gamma^\ast(S)$ is smaller or equal to one of the $\gamma(L_j)$. In particular, max$\Gamma^\ast(S)$ is
contained $\{\gamma(L_j)~|~j\}$ and hence finite. The orthants $L_j$  of rank $n$ form a complete set of representatives
of all orthants of rank $n$.\hfill$\square$

\textbf{Remark:}
We leave it to the reader to deduce that $h(S\cup S') = h(S)+h(S')$, if $S$ and $S'$ are orthohedral sets with ${\text{rk}}(S) = {\text{rk}}(S') > {\text{rk}}(S\cap S')$. 

\subsection{Piecewise isometric maps.}~ Let $S\subseteq\mathbb Z^N$ be a subset.
We call a map $f: S \to \mathbb Z^N$ \emph{piecewise Euclidean-isometric} (abbreviated as \emph{pei-map}), if $S$ is a finite disjoint 
union of orthants with the property that the restriction of $f$ to each of them is given by the restriction of an isometry. 
\\
Correspondingly, we call $f$  \emph{piecewise Euclidean-translation map} (abbreviated as \emph{pet-map}), if  $S$ is a
finite disjoint union of orthants with the property, that the restriction of $f$ to each of them is given by the
restriction of a translation.
\\
If a bijection $f: S\to S'$ is a pei-map (resp a pet-map), so is $f^{-1}$ and we say that $S$ and $S'$ are
\textit{pei-isomorphic} (resp. \textit{pet-isomorphic}). 

By the argument used in the proof of Lemma \ref{lemma2.2} above one shows that if $f$ is a pei-map, then each orthant
$L\subseteq S$ contains a commensurable suborthant on which $f$ restrict to an isometric embedding. In fact, we leave it
to the reader to observe
\begin{lemma}\label{lemma2.3} Let \mbox{$f:\mathbb Z^N \to \mathbb Z^N$} be an injective map on an orthohedral set
$S\subseteq\mathbb Z^N$. Then $f$ is a pei(resp. pet)-injection if and only if every orthant $L$ of $\mathbb Z^N$
contains a commensurable suborthant $L'\subseteq L$ on which $f$ is given by an isometry (resp translation) onto
$f(L')\subseteq\mathbb Z^N$.
\end{lemma}

It follows that every injective pei-map \mbox{$f: S \to \mathbb Z^N$} induces a rank preserving injection \mbox{$f_\ast:
\Gamma^\ast(S) \to \Gamma^\ast(f(S))$}. $f_\ast$ does not preserve the ordering -- not even if $f$ is a pet map. But
since it is rank-preserving, it does induce a bijection between the germs of maximal rank of $\Gamma^\ast(S)$ and
$\Gamma^\ast(f(S))$, whence $h(f(S)) = h(S)$.
The following observations can be left as an exercise:
\begin{lemma}\label{lemma2.4} If \mbox{$f: S \to \mathbb Z^N$} is a pet-map, then $f_\ast(\gamma)$ is parallel to $\gamma$
for each $\gamma\in \Gamma^\ast(S)$. Hence $\tau(f_\ast(\gamma)) = \tau(\gamma)$, and  $S_\tau = f(S)_\tau$. In other
words we have the commutative diagram
$$
\begin{xy}
  \xymatrix{
      \Gamma^\ast(S)  \ar[r]^{f_\ast} \ar[d]_\tau &   \Gamma^\ast(f(S)) \ar[d]_\tau \\
      \Gamma_0^\ast(S) \ar@{=}[r] & \Gamma_0^\ast(f(S)_\tau)   
  }
\end{xy}
$$\hfill$\square$
 \end{lemma}

\subsection{Normal forms.}~
Consider the disjoint union of orthants \mbox{$S =  L_1\cup L_2\cup \ldots \cup L_m$} in $\mathbb Z^N$. Assuming that
$\text{rk}S<N$ we have enough space to parallel translate each $L_i$ to an orthant $L'_i$ in such a way that the $L'_{i}$ are
still pairwise disjoint, but that each (oriented) parallelism class of the orthants $L'_i$ is assembled to a stack.
This describes a pet-bijection $S\to S' = \bigcup_jS_j$ , where the $S_j$
stand for pairwise disjoint and non-parallel stacks of orthants. We can go one step further by observing that when the maximal orthants of a stack $S_i$ are parallel to suborthants of the stack $S_j$, then there is a pet-bijection   $S_i\cup S_j \to S_j$ which feeds $S_i$ into $S_j$. Hence we can delete all stacks $S_i$ of orthants that are parallel to a suborthant of some other $S_j$ and find

\begin {proposition}\label{proposition2.5}\emph{(pet-normal form)}
Each orthohedral set $S$ is pet-isomorphic to a disjoint union of stacks of
orthants $S' =\bigcup_ j S_j$, with the property that no maximal orthant of any
$S_j$ is parallel to a suborthant in some  $S_k$, if $k\not= j.$\hfill $\square$
\end{proposition}
\begin{corollary}\label{corollary2.6}\emph{(pei-normal form)}
Each orthohedral set $S$ is pei-isomorphic to a stack of orthants.\hfill $\square$
\end{corollary}

As observed in Section 1.3, rank $\text{rk}S$ and height $h(S)$ are pei-invariant; hence they can be read of from the pei-normal
form; and the pair $(\text{rk}S, h(S))$ characterizes $S$ up to pei-isomorphism. For the corresponding pet-result we consider the
height function

\begin{tabular}{lrc} 
 & (1.1) & $ h_S: \Gamma_0^\ast \longrightarrow \mathbb N_0$
\end{tabular}

which assigns to each 0-based orthant  $L\in \Omega_0^\ast =\Gamma_0^\ast$ the number of maximal germs
$\gamma\in\text{max}\Gamma^\ast(S)$ with  $\tau(\gamma) = L$, which is finite by Lemma \ref{lemma2.2}. The support
$\text{supp}(h_S) \subseteq \Gamma_0^\ast$ is the set of all 0-based orthants $L$ with $h_S(L) > 0$. From the Exercise in Section
1.2 we infer that max$\Gamma_0^\ast(S_\tau)\subseteq \text{supp}(h_S)$, and that this is not, in general an equality. One
observes easily that the equality 

\begin{tabular}{lrc} 
 & (1.2) & $\tau(\text{max}\Gamma^\ast(S)) = \text{max}\Gamma_0^\ast(S_\tau))$ (or equivalently:
max$\Gamma_0^\ast(S_\tau) = \text{supp}(h_S)$ )
\end{tabular}

is a necessary condition for $S$ to be in pet-normal form. Thus we call $S$ \textit{quasi-normal} if the equation (1.2)
holds. Of course, a quasi-normal orthohedral set is not necessarily in pet-normal form. But as quasi-normality implies
that $\tau$ restricts to a surjection \mbox{$\tau : \text{max}\Gamma^\ast(S) \twoheadrightarrow
\text{max}\Gamma_0^\ast(S_\tau)$}, max$\Gamma^\ast(S)$ is the pairwise disjoint union of the fibers $f^{-1}(\gamma)$, which consist of  $h_S(\gamma)$  germs
parallel to $\gamma$. This can be viewed as a weak germ-version of the pet-normal form.

\begin{lemma}\label{lemma2.7}
 If \mbox{$f:\mathbb Z^N \to \mathbb Z^N$} is a pet-injection of a quasi-normal orthohedral set $S\subseteq\mathbb Z^N$,
 then $f_\ast(\emph{max}\Gamma^\ast(S)) \subseteq \emph{max}\Gamma^\ast(f(S))$.
\end{lemma}
\textit{Proof.} By Lemma \ref{lemma2.3} $f$ induces a rank preserving bijection \mbox{$f_\ast: \Gamma^\ast(S) \to
\Gamma^\ast(f(S))$}, and by Lemma \ref{lemma2.4} $f(S)_\tau = S_\tau.$ Let $\gamma\in\text{max}\Gamma^\ast(S)$. Then
we know that $\tau(\gamma)$ is maximal in $\Gamma_0^\ast(S_\tau)$. Since $f$ is a pet map, we also know that
$\tau(f_\ast(\gamma)) = \tau(\gamma)$; hence $\tau(f_\ast(\gamma))$ is maximal in  $\Gamma_0^\ast(S_\tau) =
\Gamma_0^\ast(f(S)_\tau)$. We claim that $f_\ast(\gamma)$ is maximal in $\Gamma^\ast(f(S)_\tau)$. Indeed, if $f_\ast(\gamma)$ is
not in max$\Gamma^\ast(f(S)_\tau)$, then $\tau(f_\ast(\gamma))$ cannot be maximal in $\Gamma_0^\ast(f(S)_\tau)$. This
shows that  $f_\ast(\text{max}\Gamma^\ast(S)) \subseteq \text{max}\Gamma^\ast(f(S))$, as asserted.\hfill$\square$

\begin{corollary}\label{corollar2.8}
If \mbox{$f:S \to S'$} is a pet-isomorphism between quasi-normal orthohedral sets, then $f_\ast(\emph{max}\Gamma^\ast(S))
= \emph{max}\Gamma^\ast(S')$ and $h_S = h_{S'}$.\hfill$\square$
\end{corollary}

This shows, in particular, that the stack heights in a pet-normal form are uniquely determined and characterize S up to
pet-isomorphism.

\section{Permutation groups supported on orthohedral sets}
\subsection{pei- and pet-permutation groups.}~ Let $G = \text{pei} (\mathbb Z^N)$ denote the group of all pei-permutations
of $\mathbb Z^N$.
Given any subset $S \subseteq \mathbb Z^N$, we write $G(S)$ for the subgroup of
$G$ supported on $S$, i.e. \mbox{$G(S) = \{g\in G ~|~ g(x) = x$} for every $x\in
(\mathbb Z^N - S) \}$. If $S$ is orthohedral, so is its complement and hence every pei-permutation of $S$ extends to a pei-isomorphism of  
$\mathbb Z^N$, and we write pei$(S)$ for $G(S)$ when this is convenient.

As an immediate consequence of Corollary \ref{corollary2.6} we have 
\begin{corollary}\label{corollary3.1}  If $S\subseteq\mathbb Z^N$ is an
orthohedral subset, then \emph{pei}$(S)$ is isomorphic to \emph{pei}$(S')$, where $S'$ is a stack of orthants of rank
$\emph{rk}S$ and height $h(S)$.\hfill$\square$
\end{corollary}

The set of all pet-permutations on the orthohedral set $S$ is the pet-subgroup pet$(S) \le G(S)$. As an immediate consequence of  Theorem \ref{proposition2.5} and its subsequent remark we find \begin{corollary}\label{corollary3.2} If $S\subseteq\mathbb Z^N$ is an
orthohedral subset and $S' = \bigcup_j S_j$  its pet-normal form, then
\emph{pet}$(S)$ is isomorphic to \emph{pet}$(S')$.\hfill$\square$
\end{corollary}

\subsection{The germ stabilizer.}~Let $S\subseteq\mathbb Z^N$ be an orthohedral set. We start by attaching to each germ
$\gamma \in \Gamma^\ast(S)$ the union $\langle \gamma \rangle \subseteq \mathbb Z^{N}$ of all orthants representing
$\gamma$, which we call the \textit{tangent coset at} $\gamma$. Clearly, $\langle \gamma \rangle$ isometric to $\mathbb
Z^{\text{rkL}}$. Given a pei-permutation $g\in G(S)$, we find in each orthant $L$ of $S$ a suborthant $L'$ commensurable
with $L$, such that the restriction \mbox{$g|_{L'}: L' \to S$} is an isometric embedding, and putting $g\gamma (L') =  \gamma
(gL')$ yields a well defined induced action of $G(S)$ on $\Gamma^\ast(S)$. A certain control over the stabilizer $C(\gamma )$ of $\gamma$
is given by

\begin{proposition}\label{proposition2.3}$C(\gamma )$ acts canonically on the tangent coset $\langle \gamma \rangle
\cong \mathbb Z^{\emph{rk}(\gamma)}$ by isometries. The kernel $K(\gamma ) = {\rm ker} \big(C(\gamma )\to
\emph{Isom}\langle \gamma \rangle\big)$ of this action consists of all elements $g\in G(S)$ with the property that
$\gamma$ is represented by an orthant pointwise fixed by $g$.
\end{proposition}

\textit{Proof.} Let $g\in G(S)$ with $L'\subseteq L$ as above. If $g\gamma(L') =  \gamma(gL') = \gamma(L')$, then $L'$ and
$gL'$ are commensurable and hence the restriction map $g|_{L'}$ embeds $L'$ isometrically into $\langle \gamma \rangle$.
This embedding extends to an isometry of $\langle \gamma \rangle$ and is independent of the particular choice of $L'$. One
observes that this well defines an action of $C(\gamma)$ on $\langle \gamma \rangle$; the second statement is obvious.
\hfill$\square$

\textbf{Remark.} Note that Isom$(\mathbb Z^k)$ contains the translation subgroup $\mathbb Z^k$ as a normal subgroup of
finite index.

\subsection{A normal series for $\mathbf{G(S)}$.}~ With $S$ as in 2.2 we consider the chain of subsets
	\[ (^\ast)\quad	\Gamma^\ast(S) = \Gamma[0] \supseteq  \Gamma[1] \supseteq  \Gamma[2] \supseteq   \ldots  \supseteq   \Gamma[j] \supseteq  \ldots \supseteq  \Gamma[n] = \Gamma^n(S),\] 
where $\Gamma[j] = \bigcup_{k\ge j} \Gamma^k(S)$. As the action of $G(S)$ on $\Gamma^\ast(S)$ preserves the rank it yields an ascending series for the germwise stabilizers
	\[ 1 =  C(\Gamma[0]) \subseteq  C(\Gamma[1]) \subseteq  C(\Gamma[2])  \subseteq   \ldots  \subseteq  C(\Gamma[j]) \subseteq  \ldots \subseteq  C ( \Gamma[n]) =  C \big(\Gamma^n(S) \big),\]
where
	\[ C(\Gamma[j]) := \{ g\in  G(S) ~\mid~ g\gamma  = \gamma , ~\text{for all germs}~ \gamma \in  \Gamma^\ast(S)
	~\text{with}~ \text{rk}\gamma \ge j\} .\] But we can do better: Let 
	
\begin{eqnarray}
 K(\Gamma[j] & := & \bigcap_{\gamma \in  \Gamma[j]} K(\gamma )\nonumber\\
             & =  & \{
          			g\in  C(\Gamma[j]) \mid g ~ \text{pointwise fixes in each}~ L\in \Omega^k(S),\nonumber\\
             &&   \hspace{26mm}	\text{with~} k\leq j, \text{a commensurable suborthant}. \} \nonumber
\end{eqnarray}

We claim 

\begin{proposition}\label{proposition3.4}
$K(\Gamma[j])$ acts by finitary permutations on $\Gamma[j - 1]$. 
\end{proposition}
\textit{Proof.} Let $g\in  K(\Gamma[j])$. We start by proving, that $g$ is supported on an orthohedral subset $S'\subseteq S$ with rank $S' = j-1$. We prove this by induction on $i = n - j.$ 

If $j = n$, then $g\in  K(\Gamma[n])$ asserts that each rank-$n$ germ is
represented by an orthant pointwise fixed by $g$. By Corollary \ref{lemma2.2} $\Gamma^n(S)$ is finite, and when we remove
from $S$ for each $\gamma \in  \Gamma^n(S)$ an orthant representing $\gamma$ which is pointwise fixed by $g$, we end up with a rank-$(n-1)$ orthohedral subset $S'\subseteq S$ which supports $g$. 
\\
If $j < n$ and $g\in  K(\Gamma[j])$, the inductive hypothesis asserts that g is
supported on a orthohedral subset $S'\subseteq S$ of rank $j.$ As before, $\Gamma^j(S')$ is finite. So we can
remove representing  rank-$j$ orthants pointwise fixed by $g$ and find an
orthohedral rank-$(j-1)$ subset $S''\subseteq S'$ supporting $g$.
$S''$ contains all (finitely many) rank-$(j-1)$ orthants of $S$ which are not pointwise fixed by $g$. This proves the proposition, for it is the set of all germs of such orbits that $g$ has to permute in $\Gamma[j - 1].$\hfill$\square$

The kernel of the finitary permutation representation  $K(\Gamma[j]) \to$ {\rm Sym}$(\Gamma [j - 1])$ is the subgroup
$C(\Gamma[j - 1])$; hence we have

\begin{theorem}\label{theorem2.5} For any orthohedral set $S$ with $\emph{rk}S=n$ the
group $G(S)$ admits the normal series 
\[ \begin{array}{lcl}  
1 \le K(\Gamma[1]) \le C(\Gamma[1]) \le K(\Gamma[2]) &\le& C(\Gamma[2]) \le \ldots\\
&\le& C(\Gamma[n-1]) \le K(\Gamma[n])\le  C(\Gamma[n]) \le G(S),
\end{array} \]
with the factor groups  \emph{pei}$(S)/C(\Gamma[n]) = \emph{sym}(n),$
\[ \begin{array}{rcl}  
				C(\Gamma[n])/K(\Gamma[n]) &\le& \emph{Isom} (\mathbb Z^n)^{h(S)}  ~\emph{Abelian-by-finite},\\
				 C(\Gamma[j])/K(\Gamma[j]) &\le& \emph{Isom} (\mathbb Z^j)^\infty  ~\emph{Abelian-by-(locally finite), if}~ 1<j<n\\
  			       K(\Gamma[j])/C(\Gamma[j - 1]) &\le& \emph{sym} (\infty)   ~\emph{locally finite,  if}~ 1\le j\le n.
\end{array} \]\hfill$\square$
\end{theorem}
\begin{corollary}\label{corollary 2.6} $G(S)$ is elementary-amenable. \hfill$\square$
\end{corollary}

\section{A lower bound for the finiteness length of pei$\mathbf{(S)}$}
In this section we will define a certain ``diagonal'' subgroup, ${\text{pei}}_{\rm{\text{dia}}}(S) \leq {\text{pei}}(S)$, and prove
\begin{theorem}\label{theorem4.1}  For every orthohedral set $S$ we have 
$$ fl{\rm  \big( pei}  (S) \big) \ge 
fl \big({\emph{pei}}_{\emph{dia}} (S) \big) = h(S) - 1.$$
\end{theorem}

\subsection{The height of a pei-injection ~$\mathbf{ f: S\to S.}$}~ We start with
a general observation on the set of germs, when an orthohedral set $S$ comes with
a decomposition of a disjoint union  $S = A\cup B$ of two orthohedral subsets.
In that case every orthant $L\subseteq S$ inherits the decomposition  $L = 
(A\cap L)\cup (B\cap L)$, which shows that one of the orthants of either $A$ or $B$ is commensurable to $L$. 
This shows that the germs of S have an induced disjoint decomposition $\Gamma^k(S)  =  \Gamma^k(A)\cup \Gamma^k(B)$ for
each  $k.$

Now let $S$ be an orthohedral set of rank ${\text{rk}}S = n.$ We can represent the
rank-$n$ germs of $S$ by pairwise disjoint orthants $L_1, \ldots , L_h, h =
h(S)$, with the property that the restriction of $f$ to each $L_i$ is an
isometric embedding into $S$. $f(L_i)$ is then commensurable to some $Lj$, and
since $f$ is injective it follows:
\\ $f$ permutes the germs  $\gamma (L_1),\ldots , \gamma (L_p)$, and $ {\text{rk}}(S - f(S)) < {\text{rk}}S.$
\\
As $S - f(S)$ is an orthohedral set, we now obtain  that the number of
rank-$(n-1)$ germs in $S - f(S)$ is finite. We call this number \textit{the
height of} $f$, denoted $h(f) = h\big( S - f(S)\big) = h(S - Sf).$

\begin{lemma}\label{lemma3.2}
\begin{enumerate} \item If $f, g: S\to S$ are two pei-injections, then $h(f\cdot g) = h(f) + h(g).$
\item If  $A\subseteq S$ is an orthohedral subset whose complement $A^c = S - A
$ has rank ${\emph{rk}}A^c < n = {\emph{rk}}S$, then the height of any pei-injection $f: S\to S$
is given by  $h(f) = $ \mbox{$h\big( A\cap f(A)^c \big)$} $-$ $h\big( A^c\cap
f(A)\big).$
\end{enumerate}
\end{lemma}

\textit{Proof.}
i) Consider the disjoint union $S = (S - Sg)\cup Sg$. As $f$ is injective $Sf
=$ \mbox{$(Sf - Sgf)$} $\cup~gf$ is also a disjoint union. Hence so is  $S = (S - Sf)\cup Sf = (S - Sf)\cup$ \mbox{$(Sf -Sgf)\cup Sgf$}, 
and we find \[	S - Sgf = (S - Sf)\cup (Sf - Sgf). \] Now, $f$ is a pei-bijection
between $(S - Sg)$ and $(S - Sg)f =  (Sf - Sgf)$; and a pei-bijection of a an orthohedral set induces a pei-bijection on its germs. Thus the number of rank-$(n-1)$ germs of $(S - Sg)$ and $(Sf - Sgf)$ are the same. This proves i).
\\
ii) Each pei-injection $f: S\to S$ induces an injection $f^\ast:
\Gamma^{n-1}(S)\to \Gamma^{n-1}(S)$. We abbreviate $B = A^c$ and know from ${\text{rk}}B < n$  that $\Gamma^{n-1}(B)$ is finite. Hence $f^\ast$ restricts to a bijection  $f^\ast: \Gamma^{n-1}(B)\to \Gamma^{n-1}\big(f(B)\big)$. On the complement we find the induced injection $f^\ast: \Gamma^{n-1}(A) \to  \Gamma^{n-1}\big(f(A)\big).$
\\We use the abbreviation $P^\ast:=\Gamma^{n-1}(P) $ for $P = S, A, B$, and
consider the disjoint union \[ \begin{array}{rcl}	
	S^\ast - f^\ast(S^\ast) &	=&\!\!\!\! \big(A^\ast - A^\ast\cap f^\ast(S^\ast)\big) \cup  \big(B^\ast - B^\ast\cap f^\ast(S^\ast)\big)\\
		&	=  & \!\!\!\! \big( A^\ast - A^\ast\cap f^\ast(A^\ast) - A^\ast\cap f^\ast(B^\ast) \big) \cup \big( B^\ast - B^\ast\cap f^\ast(A^\ast) - B^\ast\cap f^\ast(B^\ast) \big).
\end{array}\]
Using the fact that $B^\ast$ and $A^\ast - f^\ast(A^\ast)$ are finite, we find
for $h(f) = $ \mbox{$h\big( S^\ast - f^\ast(S^\ast)\big)$} \[  h(f)  =  h \big(
A^\ast- A^\ast\cap f^\ast(A^\ast)\big) - h \big( A^\ast\cap f^\ast(B^\ast) \big) + h(B^\ast) - h \big( B^\ast\cap f^\ast(A^\ast)\big) - h \big( B^\ast\cap f^\ast(B^\ast)\big) .\] Now we apply that $h(B^\ast) = h\big( f(B^\ast) \big)$ and observe that
\[ 
\begin{array}{lcl}
- h \big( A^\ast \cap f^\ast(B^\ast)\big) &+& h(B^\ast) - h \big( B^\ast\cap f^\ast(B^\ast)\big)\\
 &=& - h\big( A^\ast\cap f^\ast(B^\ast)\big) + h\big( f(B^\ast)\big) - h \big( B^\ast\cap f^\ast(B^\ast)\big)\\
 &=& h\big( f(B^\ast) - h\big((A^\ast\cap f^\ast(B^\ast) \big)\cup\big( B^\ast\cap f^\ast(B^\ast)\big)\big) = 0.
\end{array}
\]

Hence our expression for h(f) simplifies to
\[	h(f)  =  h \big( A^\ast- A^\ast\cap f^\ast(A^\ast) \big)  - h \big( B^\ast\cap f^\ast(A^\ast) \big) =  h \big( A^\ast\cap f^\ast(A^\ast)c\big)  - h \big(B^\ast\cap f^\ast(A^\ast)\big)\]
as asserted. \hfill$\square$
\vskip 2mm 
\subsection{Monoids of pei-injections.}~
Let $S$ be an orthohedral set of rank $n = {\text{rk}}S$ in pet-normal form.  In
particular $S$ is the pairwise disjoint union of finitely many specified stacks of orthants. Let 
max$\Gamma^\ast(S)$ denote the finite set of all maximal germs of S.
We write  $M(S)$ for the monoid of all pei-injections $S \to S$. It is endowed with the height function  $h: M(S)\to
\mathbb N$ of Section 3.1. Let $M_0(S)$ be the submonoid of all pei-endoinjections of $S$, which fix all maximal germs of $S.~ M_0(S)$
is of  finite index in $M(S)$ since max$\Gamma ^\ast(S)$ is finite.  As in the proof of Proposition 2.2 we see, that each 
$f\in M_0$ induces an isometry $\tau_{(f,\gamma)}: \langle \gamma \rangle \to \langle \gamma \rangle$ on the tangent
coset of each germ $\gamma\in \max\Gamma^\ast(S).$ Thus we  have a homomorphism

\begin{tabular}{rc} 
(3.1) & $ \quad	\kappa : M_0(S)\longrightarrow >\bigoplus_{\gamma\in \max\Gamma^\ast(S)} {\rm Iso} \left(\langle \gamma\rangle \right), ~ {\rm given~ by}~ \kappa (f) =
\bigoplus_{\gamma\in \max\Gamma^\ast(S)} \tau_{(f,\gamma)}$.
\end{tabular}

The \textit{translation submonoid} $M_{tr}(S)\subseteq M_0(S)$ consists of all
$f\in M_0(S)$ with the property that the induced maps $\tau_{(f,\gamma)}: \langle \gamma \rangle \to \langle \gamma \rangle$
are translations for each $\gamma\in \max\Gamma^\ast(S)$.
Since the translation subgroup of  Isom$(\langle \gamma \rangle)$ is of finite index, $M_{tr}(S)$ has finite index in
$M_0(S)$. And restricting (3.1) yields a surjective homomorphism
 
\begin{tabular}{rc}
(3.2) & $\quad	\kappa : M_{\rm tr}(S)\longrightarrow >\bigoplus_{\gamma\in \max\Gamma^\ast(S)} \mathbb
Z^{\rm \text{rk}\gamma} = \mathbb Z^N, {\rm with}~ N = \Sigma_{\gamma\in \max\Gamma^\ast(S)} \text{rk}\gamma.$
\end{tabular}

Every orthant $L$ contains a characteristic \textit{diagonal element} $u_L\in
L$: the sum of the canonical basis of $L$. We write $t_L: L\to L$ for the
translation given by addition of $u_L$ and call this the \textit{diagonal
unit-translation} of $L$. The general diagonal translations on $L$ (i.e on $\langle L \rangle$) are given by addition of
an integral multiple  of $u_L$. By the \textit{diagonal submonoid}  $M_{\text{dia}}(S)\subseteq M_{\rm tr}(S)$ we mean the set of all
elements $f\in M_{\rm tr}(S)$ with the property, that for each $\gamma \in \max\Gamma^\ast(S)$ the induced isometry $\tau_{(f,\gamma)}:
\langle \gamma \rangle \to \langle \gamma \rangle$ is a diagonal translation. Restricting (3.2) yields the homomorphism

\begin{tabular}{rc}
 (3.3) &  $\quad	\kappa : M_{\text{dia}}(S)\longrightarrow >\bigoplus_{\gamma\in \max\Gamma^\ast(S)} \mathbb Z
 = \mathbb Z^{\mid \max\Gamma^\ast(S) \mid}$.
 \end{tabular}
 
 We write max$\Omega^*(S)$ for the set of all maximal orthants of the stacks of S, and consider the set $T = \{t_L\mid
 L\in \max\Omega^\ast(S)\}$ of all diagonal unit-translations of these orthants. Each $t_L \in T$ extends canonically to a pei-injection on \mbox{$t_L : S
 \to S$}, which is the identity on $S - L$.
 We denote it by the same symbol $t_L $, and with this interpretation $T$ generates a free-Abelian submonoid mon$(T)\le M_{\text{dia}} (S)$.
\begin{definition}
Given $f, f '\in M_{\text{dia}} (S)$ we define $f\le f '$ if there is some $t\in {\rm mon}(T)$ with $tf = f '$.
\end{definition}
\textbf{Observation.} $M_{\text{dia}} (S)$ is a \emph{directed partially ordered set}.\hfill$\square$

It is an important fact, that the height function $h: M(S) \to  \mathbb N$  is
order preserving and its restrictions to totally ordered subsets of $M(S)$ are
injective.  We will also have to consider slices of  $M_{\text{dia}} (S)$. For
given $r,s \in \mathbb N_0$, $r \le s$ we put 
\begin{eqnarray} 
M^{[r,s]} &:=& \{f\in M_{\text{dia}} (S)\mid r \le h(f) \le s\} ~~~ \text{and}\nonumber\\
M^{[r,\infty]} &:=& \{f\in M_{\text{dia}} (S)\mid r \le h(f)\}.\nonumber
\end{eqnarray}
 $M^{[r,\infty]}$ inherits the partial ordering from $M_{\text{dia}} (S)$ and is
also a directed set.
\vskip 2mm 
\subsection{Maximal elements $\mathbf{< f}$ in $\mathbf{M_{\text{dia}} (S)}$.}~
From now on we assume that all maximal orthants of the stacks of $S$ have the same finite rank
$n = \text{rk}S$. We put $\Lambda := \max\Omega^\ast(S)$. We write
\[	M_{< f} = \{a\in M_{\text{dia}} (S)\mid a< f\}, \quad M_{\le f} = \{a\in M_{\text{dia}} (S)\mid a\le f\}  \] 
for the \textit{``open resp. closed cones below $f$''} and aim to understand the set of all maximal elements
of $M_{< f}$ .
For this it will be convenient to have to introduce the abbreviation for the points on the (finite)
boundary of the maximal orthants $L$, so we set \mbox{$\partial L := L - Lt_L$}.

\begin{lemma}\label{lemma3.4}
Let $b$ be a maximal element of  $M_{< f}$. Then there is a unique maximal
orthant $L \in \Lambda$ with the property that $f = t_L b$ and $h(f) = h(b) + n$.
Furthermore $b$ is given as the union $b = b'\cup b''$, where $b': \partial L \to (S - Sf)$ is a pei-injection, and $b'':
(S - \partial L) \to Sf$ is the restriction $(t_L^{ -1}f)\mid_{ (S - \partial L)}$.
Converseley, if $c': \partial L \to (S -  Sf)$ is an arbitrary pei-injection distinct to
$b'$, then the union $c = c'\cup b''$ is a maximal element of $M_{< f}$
distinct to $b$.
\end{lemma}

\textit{Proof.} For each element $b\in M_{< f}$ there is some $t\in {\rm mon}(T)$ with $f = tb.~ t$ has a unique reduced expansion as a product of elements 
of $T$; let $l(t)$ denote the length of this expansion. As $h(t_L ) = L = n$ for each $L\in \Lambda$  we have $h(t) = n~l(t)$. It follows that if 
$b$ is maximal, then $h(t) = n$ and $t = t_L \in T$ for some $L\in \Lambda$.
The maximal orthant $L$ is uniquely determined by the fact that the restriction of $f$ and $b$ coincide on $S - L.$ 
The restriction $b''$ of $b$ to $(S - \partial L)$ coincides with $(t_L^{
-1}f)\mid_{ (S - \partial L)}$, and has its image in Sf. The restriction $b'$ of $b$ to  $\partial L$ is not determined by $f$ and $L.$ As $b$ and $f$ are injective we know that \[ \begin{array}{rcl} \emptyset  = (\partial L)b\cap (S - \partial L)b &=& (\partial L)b\cap \big( (S - L)b\cup Lt_L b\big)\\
				& =& (\partial L)b\cap \big( (S - L)f\cup Lf\big)\\ 
			&	 = &(\partial L)b\cap Sf.
\end{array} \]

Hence $b'$ can be viewed as a pei-injection $\partial L\to (S - Sf)$.
If we replace $b'$ by another pei-map $c': \partial L\to (S - Sf)$, the union
$c = c'\cup b''$ will still satisfy  $f = t_L c$ and $h(f) = h(c) + n$. This shows that $c$ will also be maximal in $M_{< f} .$\hfill$\square$

\begin{lemma}\label{lemma3.5}
Let $B\subseteq M_{< f}$ be a finite set of maximal elements of $M_{< f}$. Then the following conditions are equivalent:
\begin{enumerate}
\item   The elements of $B$ have a common lower bound $\delta$ in $M_{< f}$
\item For every pair $(b,b')\in B\times B$ with $b\not= b'$ and $tb=f =t'b'$ for diagonal unit-translations $t$, $t'$, we
have
\begin{enumerate}
\item  $t\not= t'$  and 
\item 	  $b(\partial L)\cap b'(\partial L') = \emptyset$, where $L$ resp. $L'$ are the maximal orthants of $S$ on which
$t$ resp. $t'$ acts non-trivially.
\end{enumerate}
\end{enumerate}
\end{lemma}

\textit{Proof.} i) $\Rightarrow$ ii). Let $\delta$ be a common lower bound of the elements
of $B$. Then for every pair $(b,b')\in B\times B$ there are diagonal translations $d, d'\in {\text{mon}}(T)$ with $d\delta
= b$ and $d'\delta = b'$.
From $tb = f = t'b'$ we obtain  $td\delta =  t'd'\delta$ and  conclude $td = t'd'$.
The assumption $t = t'$ would now imply $d = d'$ and hence $b = b'$.
\\
Let $L$ resp. $L'$ denote the maximal orthants of $S$ on which $t$ resp. $t'$ acts non-trivially. As $d, d'$ are diagonal translations, we have 
$d(L)\subseteq L$ and $d'(L')\subseteq L'$. From $t \not=  t'$ we know $L\cap L' = \emptyset$ , and hence $(\partial L)d\cap (\partial L')d' = \emptyset$. 
Since $\delta$ is injective, this implies  $\emptyset = (\partial L)d\delta\cap (\partial L')d'\delta = (\partial L)b\cap (\partial L')b'$, as asserted.

ii)$\Leftarrow$ i). For each $b\in B$ we have some diagonal unit-translation $t_b\in T$ with $t_b b = f$, and we put

\begin{tabular}{r c}
(3.4)& $t_B := \prod_{ b\in B}t_b$.
\end{tabular}

By assumption (i) the maximal orthants $L_b $ on which $t_b $ is a diagonal unit-translation are pairwise disjoint. 
Thus $|B|$ $\le h(S)$, and $S$ decomposes in the disjoint union  $S  =
\left(\bigcup_{b\in B}L_b \right) \cup  S'$.
We define the pei-injection $\delta_B: S\to S$ as follows:
$$ \delta_B := \begin{cases}
{t_b}^{-1} f & \text{on each } L_b t_b \\ 
b & \text{on the complements } \partial L_b = L_b - L_b t_b \\ 
f & \text{on } S'. 
\end{cases} $$

Assumption (ii) guarantees that the restriction of $\delta_B $ to the union 
$\bigcup_{ b\in B}\partial L_b  = \bigcup_{ b\in B} (L_b  - L_b t_b ) = (S - St_B)$ is injective.  And since the image of each $\partial L_b$ is disjoint to
$f(S)$, we also find that the image of $(S - St_B)$ is disjoint to $f(S)$, and also to $f(S')\subseteq f(S)$. This
shows that $\delta$  is a pei-injection. It remains to prove that $\delta_B $ is a common lower bound for the elements of $B$.
By commutativity we find elements $s_b\in ~{\rm mon}~(T)$ with $t_B \delta_B  = t_b s_b\delta_B $, where $s_b = \prod_ { x\in (B - \{b\})}t_x$.
One observes that $s_b\delta_B $ and $b$ agree on $(S- St_b) = (L_b  - L_b t_b)$, and that $t_B\delta_B = f = t_b b$. 
Hence $b$ and $s_b\delta_B $ agree on $S$. 
\hfill$\square$

\begin{lemma}\label{lemma3.6}
In the situation of Lemma \ref{lemma3.5} we have for the lower bound $\delta_B$ defined in the proof:
\begin{enumerate}
\item $\delta_B$  is, in fact, a largest common lower bound of the elements of $B$
\item $h(\delta_B ) \ge h(f) - h(S)n$.
\end{enumerate}
\end{lemma}
\textit{Proof.} i) We compare an arbitrary common lower bound $\gamma$ with
$\delta_B$, the lower bound constructed in the proof above. Thus for each $b\in
B$ we are given $u_b\in {\rm mon} (T)$ with $u_b \gamma  = b.$ We fix a base element $b'\in B$ and define the diagonal translation 
$t'\in {\rm mon} (T)$ by its action on $S$ as $$ t' := \begin{cases} u_{b'} & \text{on } S' \\ 
u_b & \text{on each } L_b.
\end{cases} $$
We use the elements $s_b$ defined in the proof of Lemma \ref{lemma3.5} and observes that

\begin{tabular}{r l l}
	& $xt'\gamma  = xu_{b'} \gamma  = xb' = xs_{b'}\delta_B  = x\delta_B$ & for
	$x\in S'$\\
	& $xt'\gamma  = xu_b \gamma  = xb = xs_b\delta_B  = x\delta_B$ & for $x\in L_b$.
\end{tabular}

This shows that  $t'\gamma  = \delta_B$ , hence  $\delta_B  \ge \gamma$.
\\
ii) For the translation $t$ of (3.4) we have $t\delta_B  = f$  and can deduce
that $h(\delta_B ) =$ \mbox{$h(f) -  h(t)$} $=$ \mbox{$h(f) -  |B| n$} $\geq$ $h(f) - h(S)n.$\hfill
$\square$
\vskip 2mm 
\subsection{The simplicial complex of $\mathbf{M_\mathbf{\text{dia}} \mathbf{(S)}}$.}~
We consider the simplicial complex $|M_{\text{dia}}(S)|$,  whose vertices are the elements of $M_{\text{dia}} (S)$ and whose chains of length $k,  
a_0 < a_1 < \ldots < a_k$, are the $k$-simplices. As the partial ordering on
$M_{\text{dia}} (S)$ is directed, $|M_{\text{dia}}(S)|$ is contractible.

In this section we aim to prove

\begin{lemma}\label{lemma3.7}
If $h(f)\ge 2\cdot \emph{rk}S\cdot h(S)$, then $|M_{< f}|$  has the homotopy type of a bouquet of $\big( h(S) - 1
\big)$-spheres.
\end{lemma}

The first step towards proving Lemma \ref{lemma3.7} is to consider the covering of  $|M_{< f}|$  by the subcomplexes $|M_{\le b}|$ , where $b$ runs 
through the maximal elements of  $M_{< f}$. We write $N(f)$ for the nerve of this covering. Lemma \ref{lemma3.5} and
\ref{lemma3.6} show, that all finite intersections of such subcomplexes $|M_{\le b}|$ are again cones and hence
contractible. It is a well known fact that in this situation the space is homotopy equivalent to the nerve of the covering. Hence we have \[ |M_{<f}|  ~\textit{is homotopy equivalent to the nerve}~ N(f) \] and it remains
to compute the homotopy type of $N(f)$.

The next step is to use the results of Section 3.3 to find a combinatorial model
for the nerve $N(f)$. The vertices of $N(f)$ are the maximal elements of $M_{<f}
$, and hence, by Lemma \ref{lemma3.4}, in 1-1-correspondence to the disjoint
union  $A = \bigcup_{L\in \Lambda} A_L $, where $A_L $ stands for the set of all pei-injections  $a: \partial L\to (S - Sf ).$ 
Lemma \ref{lemma3.4} allows to translate the simplicial structure of $N(f)$ into a
simplicial complex $\Sigma(f)$ on $A$:
\\
the $p$-simplices of $N(f)$ are the $p$-element sets of maximal elements $B\subseteq M_{<f} $ with a common lower bound, and the corresponding
p-simplices of $\Sigma(f)$ are the sequences $(a_L)_{L\in \Lambda '}$, where $\Lambda '$ is a $p$-element subset of $\Lambda$  and $a_L\in A_L $ with 
the property \[ (\ast)  ~~	~\textit{The intersections of the images}~ a_L(\partial L), L\in \Lambda ', ~\textit{are pairwise disjoint.}\]

Following \cite{br87}  and \cite{sa92} we will now be able to compute
the homotopy type of $\Sigma(f)$ by means of a lemma on colored graphs.

Let $\Gamma  = (V,E)$ be a combinatorial graph, given by a set $V$ of \textit{vertices} and set $E$ of \textit{edges},
where an edge is a set consisting of two non-equal vertices. A \textit{clique} of $\Gamma$  is any subset $C\subseteq V$ with the property, that any
two vertices of $C$ are joined by an edge of $\Gamma$. The \textit{flag-complex} $K(\Gamma )$ is the simplicial complex on $V$ whose
$p$-simplices are the cliques consisting of $p+1$ vertices of $V$.  Our main example here is the complex $\Sigma(f)$, which is easily seen to be the 
flag-complex of its $1$-skeleton $\Gamma (f)$.

Let $h$ be a natural number. By an $h$-colored graph $\Gamma_h$  whose vertex set $V$ is the disjoint union of n subsets $V_1, 
\ldots V_h$. Vertices  $v, v' of V$ are said to have the same color $i$ if they are contained in the same subset $V_i$. And we 
assume that two vertices joined by an edge are never of the same color.

\begin{lemma}\label{lemma3.8}
If all colors $i\in \{1, 2, \ldots , h\}$ of an $h$-colored graph $\Gamma_h = (V,E)$ satisfy the two properties
\begin{enumerate}
\item[\emph{(1)}] $V_i$ contains at least two elements, and
\item[\emph{(2)}]  For any choice of $2(h - 1)$ vertices $u^1, \ldots u^{2(h - 1)}$ in $V - V_i $ there are two vertices $v, w$ in 
$V_i $ which are adjacent to each $u^j$, i.e., for each $j\in \{1, 2, \ldots 2(h - 1)\}$ there is an edge path of length
$2$ in $\Gamma_h$ joining $v$ and $w$ via $u^j$.
\end{enumerate}
Then the flag-complex $K(\Gamma_h)$ has the homotopy type of a bouquet of $(h - 1)$-spheres.
\end{lemma}

\textit{Proof.} We use induction on $h$, starting with the observation that the statement is trivial when $h=1$. For $h\ge 2$ we assume that $K(\Gamma_{ h-1})$ is homotopy equivalent to a bouquet of  $(h - 2)$-spheres, if $\Gamma_{h-1}$ is an $(h - 1)$-colored graph which satisfies the properties (1) and (2). We construct $K(\Gamma_h)$  in several steps, similar to the method applied in Brown's proof for Hougthon's Groups \cite{br87}. We start with choosing a base vertex $v_1\in V1$ and consider its star in $K(\Gamma_h)$, 
		\[ K_0  :=  \text{st}_{K(\Gamma h)}(v_1)\]
Then we proceed with $i = 1, 2, \ldots , h$ by taking the union of $K_{i-1} \cup V_i '$, where $V_i '$ is the set of all vertices of $V_i$  which are not joined with the base vertex $v_1$ by an edge. And we put  
			\[ K_i := ~\text{full subcomplex of}~ K(\Gamma_h) ~\text{generated by}~ K_{i-1} \cup V_i '.\]

One observes that $(V_1 \cup  \ldots \cup V_i )\subseteq K_i $ and  $(V_{i +1}\cup  \ldots \cup V_h)\cap K_i  = K_0$. In particular, $K_h = K(\Gamma_h).~ K_i $ is obtained from $K_{i -1}$ by adjoining vertices $v\in V_i$   that are not connected to the base vertex $v_1$ by an edge; then taking the full subcomplex of $K(\Gamma_h)$. Thus $K_i $ is obtained from $K_{i -1}$ by adjoining for these vertices $v$ the cone over
			\[ \text{lk}(K_{i -1} ,v) := ~\text{the link of}~ v ~\text{in}~ K_{i -1}. \]

The $1$-skeleton of $\text{lk} (K_{i -1} ,v)$ has vertex set
\[ \begin{array}{rcl}
		W &=& W_1\cup  \ldots \cup W_{i-1}\cup W_{i+1}\cup  \ldots \cup W_h ~\text{   with}\\ 
		W_j &:=& \text{set of vertices of}~ V_j  ~\text{which are joined with}~ v~\text{ by an edge},\\
&& \text{for}~  j =  1, \ldots , i-1\\
		W_j &:=& \text{set of vertices of}~ V_j  ~\text{ which are joined with}~ v~\text{ and}~ v_1~\text{ by an}\\
&& \text{edge, for}~j =  i+1, \ldots , h.
\end{array}\]

Thus, the $1$-skeleton of $\text{lk}(K_{i-1} ,v)$  is an $(h - 1)$-colored subgraph $\Gamma_{h-1}$ of $\Gamma_h$ with vertex set $W$ and colors 
$\{ 1, 2, \ldots , h\} -\{ i\} $, and $\text{lk}(K_{i-1} ,v)$  is the flag-complex $K(\Gamma_{h-1})$. Now we consider any $2(h -2)$ vertices 
$u^1 , \ldots, u^{2(h - 2)}$  of  $W - W_j$   with colors in $\{ 1, 2, \ldots , h\} -\{ i, j\} $ for some $j \in \{ 1, 2, \ldots , h\} -\{ i\} $. 
Together with the vertices $v_1$ and $v$, we obtain $2(h - 1)$ vertices $u^1 , \ldots , u^{2(h - 2)}, v_1, v$ of  $V - V_j$.  By property (2) of 
$\Gamma_h$, there exists two vertices $w, w'$ in $V_j$ , which can be joined by an edge path of length $2$ via $u^k$ for each $k \in \{ 1, 2, \ldots , 2(h - 2)\} $, 
and additionally via $v_1, v$. In particular, $w$ and $w'$ can be joined by an edge with $v_1$ and $v$, and so they are vertices of  $W_j$.  
Hence $\Gamma_{h-1}$ satisfies the two properties of the lemma, and in view of the inductive hypothesis, $\text{lk}(K_{i-1} ,v)$  is homotopy equivalent to a 
bouquet of  $(h - 2)$-spheres.

From here we can use the same arguments as in the proof of Lemma 5.3 in \cite{br87}: Starting with the contractible complex $K_0,  K_1$  is obtained
from $K_0$ by adjoining for each vertex  $v\in V_1 '$ a cone over $\text{lk}(K_0 ,v)$.  Using the homotopy type of $\text{lk}(K_0 ,v)$, we can deduce that $K_1$ is
homotopy equivalent to a bouquet of $(h - 1)$-spheres. For the next steps in the construction of $K_h$, we know that $K_i$
is obtained from $K_{i-1}$ by adjoining for each vertex $v\in V_i '$ a cone over $\text{lk}(K_{i-1} ,v)$.  In view of the homotopy type of $\text{lk}(K_{i-1} ,v)$,  we see that, up
to homotopy, the passage from $K_{i-1}$ to $K_i$   consists of the adjunction of $(h - 1)$-cells to a bouquet of $(h - 1)$-spheres. \hfill$\square $

We will now apply Lemma \ref{lemma3.8} to the $1$-skeleton  $\Gamma (f)$ of $\Sigma(f)$.
By definition its vertex set is the disjoint union $A = \bigcup_{ L\in \Lambda}
A_L$, and we regard the various $A_L$  as the coloring of  $\Gamma (f)$. The edges of $\Gamma$  are the pairs of such pei-injections $\{ a_L, a'_{L'}\}$  with disjoint images.
Thus $\Gamma (f)$ is an $h(S)$-colored graph $\Gamma (f)_{h(S)}$ in the sense above, and in order to establish Lemma \ref{lemma3.7} it remains to prove

\begin{lemma}
If  $h(f) \ge 2\cdot \emph{rk}S \cdot h(S)$, then $\Gamma (f)_{h(S)}$ satisfies the assumptions of Lemma \ref{lemma3.8}.
\end{lemma}
\textit{Proof.} Let $n := S$ and $h := h(S)$. Assumption (1) is a consequence of assumption (2) except in the trivial case
$h = 1.$
 
To prove (2) we fix $L\in \Lambda$  and consider a set of $2(h - 1)$ elements $F\subseteq A - A_L$ . We have to show that
there are two elements $a, b\in A_L$  with the property that for each $c\in F$ the image im$(c) = c\left( \partial L(c)\right)$ is disjoint to  both $a(\partial L)$ and
$b(\partial L)$. In other words: there are two pei-injections 
\[   a, b:  \partial L  \to  \big(S - Sf \big) - \Big(\bigcup_{ c\in F}{\rm im}(c)\Big). \] 
To show this it suffices to compare the height function - i.e., 
the number of rank-$(n-1)$ germs -  of domain and target. Clearly, $h\left(a(\partial L))\right) = h(\partial L) = n$, and the same applies to every 
vertex of $A$. Hence $h\left(\bigcup_{ c\in F}{\rm im}(c)\right) \le 2(h - 1)n.$ By assumption $h\big(S - Sf \big) \ge 2hn$, and so the target orthohedral 
set has height at least $2hn - 2(h - 1)n = 2 n$, which is more than the height $h(\partial L) = n$ of the domain. In this  situation one observes easily
that there are arbitrarily many different pei-injections in $A_L$ whose image is
disjoint to  $\bigcup_{ c\in F}{\rm im}(c)$. This proves the lemma. \hfill$\square$

\begin{remark}
\emph{
If we replace $M_{<f}$ by the subset $M_{r,f} :=\{a \in M_\text{dia}(S) \mid h(a) \geq r ~ \text{and} ~ a < f \}$, the
assertion of Lemma 3.7 holds true, provided $f$ satisfies the additional condition $h(f) \geq r + h(S)$.
In this case we know by Lemma 3.6 that $h(\delta_B) \geq r$, where $\delta_B$ stands for the largest lower bound of a
finite set $B$ of maximal elements of $M_{r,f}$. Thus $\delta_B$ is an element of $M_{r,f}$ and the proof of Lemma 3.7
works the same way for the reduced simplicial complex $|M_{r,f}|$.
}
\end{remark}

\subsection{Stabilizers and cocompact skeletons of $\mathbf{ M(S) }$.}~
 The group $G(S)$ of all pei-permutations acts on $M(S)$ from the right, and as $h(g) = 0$ for all $g\in G(S)$ the height 
 function $h: M(S) \to  \mathbb N$  is invariant under this action. Correspondingly, $G_{\#}(S) := G(S)\cap M_{\#} (S)$ acts on 
 $M_{\#} (S)$, where $\# $ stands for $0,$ tr, or dia.  We will also restrict attention to the various $G_{\#} (S)$-invariant 
 subsets  $M^{[r,s]}_{\#}(S)= \{ f\in M_{\#} (S)\mid r\le h(f)\le s\} $ for prescribed numbers $r \le s$ in $\mathbb
 N_0$. And also, mutatis mutandis, for the corresponding pet-groups pet$_{\#}(S)$ - note that $\text{pet}_0(S) =
 \text{pet}_\text{tr}(S)$.

We start with the following easy observation:
\begin{lemma} Two elements $f, f '\in  M(S)$ are in the same \emph{pei}$(S)$-orbit if and only if  $(S - Sf )$ and $(S -
Sf ')$ are pei-isomorphic.
\end{lemma}
\textit{Proof.} As both $Sf$ and $Sf '$ are pei-isomorphic to $S$ there is a pei-isomorphism \mbox{$g': Sf \to  Sf '$}.
If there is also a pei-isomorphism $g'': (S - Sf ) \to  (S - Sf ')$ then the union $g = g'\cup g''$ is a pei-permutation of 
$S$ with $fg = f '$. Conversely, $fg = f '$ implies $(S - Sf ') = (Sg - Sfg) = (S - Sf)g$, hence $(S - Sf ')$ is pei-isomorphic 
to $(S - Sf). \hfill\square$

Since orthohedral sets of the same rank and height are pei-isomorphic by Corollary 1.5, it follows that pei$(S)$ acts
transitively on the set of all pei-injections of a given rk$(S - Sf)$ and height $k$. The very same can be said for the 
action of  $G_{\#} (S)$ on $M_{\#} (S).$

Let  $\Delta  = (a_0 < a_1 < \ldots < a_{k-1} < a_k)$ be a  $k$-simplex of $ |M(S)|$. By definition there are elements 
$t_1, t_2, \ldots , t_k \in {\rm mon} (T)$, with $a_i = t_i a_0$ for all $i$; they are uniquely defined and form a $k$-simplex   
$ \Delta ' = (id<t_1< \ldots <t_{k-1}<t_k) \in |{\rm mon}(T)|$ . Moreover, putting $\sigma (\Delta ) := (\Delta , a_0)$
defines a bijection 
			\[\sigma :~   |M(S)| ~ \longrightarrow ~  |{\rm mon}(T)| \times ~ M(S).\]

The action of pei$(S)$ on $ |M(S)|$  is given by $(a_0 <a_1 < \ldots <a_k )g = (a_0 g< a_1 g< \ldots < a_{k-1} g< a_k g)$.
We can leave it to the reader to observe that this action induces, via $\sigma$, on $ |{\rm mon} (T)| \times M(S)$ the
$G(S)$-action given by simple right action on $M(S).$

The simple structure of the $G_{\#} (S)$-action on $ |M_{\#} (S)|$  has two immediate consequences:

\begin{corollary} \begin{enumerate} \item The stabilizer of a $k$-simplex of $ |M_{\#} (S)|$  coincides with the
stabilizer of its minimal vertex $f$ and is isomorphic to $G_{\#} (S - Sf).$
\item  For every numbers $r \le s$ in $\mathbb N_0$ the simplicial complex of  $M^{[r,s]}_{\#}(S)= \{ f\in
M_{\#} (S)\mid r\le h(f)\le s\}$ is cocompact under the $G_{\#} (S)$-action
\end{enumerate} \end{corollary}

\textit{Proof.}
i) One observes that right action of $g\in G_{\#} (S)$ on $M_{\#} (S)$ fixes an element $f\in M_{\#} (S)$
if and only if $g$ restricted to $Sf$ is the identity. In other words, the stabilizer of the vertex $f\in M_{\#} (S)$ is isomorphic to $G_{\#} (S - Sf).$
\\
ii) We use the interpretation of a simplex $\Delta  = (a_0 < a_1 < \ldots < a_{k-1} < a_k) \in  |M(S)| $ in $|{\rm mon} (T)| \times~ M(S)$. 
Since $G_{\#} (S)$ acts transitively on the set of all pei-injections in $M_{\#} (S)$ of a given
rk$(S-Sf)$ and height $k$, the bound on $h(a_0)$ allows only finitely many $G_{\#} (S)$-orbits on the second component
$M(S)$. The bound on $h(a_i)$ for $i=1,\ldots,k$ allows only finitely many simplices in the first component \mbox{$| {\rm mon} (T)|$}.
\hfill$\square$
\vskip 2mm 
\subsection{The conclusion.}~ Here we put things together to prove Theorem 3.1, i.e
$fl \big(G(S)\big) \ge fl\big(G_{\text{dia}}(S)\big) = h(S) - 1.$

\textit{Proof.} We will first show, by induction on $n = S$, that $fl\big(G_{\text{dia}}(S)\big) = h(S) - 1.$ If  $n = 1$,
then the group $G_0(S)$ is the Houghton group on $h(S)$ rays and has finite index in $G(S).$  In that case the assertion
is due  to Ken Brown \cite{br87}.

Now we assume $n > 1$. Here we use $M^{[r,s]} = \{ f\in M_{\text{dia}} (S)\mid r\le h(f)\le s\} , r,s\in \mathbb N$. 
Since $f\in M^{[r,s]}$ is a diagonal pei-injection, the height of $f$ is a multiple of $n$.
So we fix the lower bound $r = nk_0$, $k_0 \in \mathbb N$, and consider the filtration of  $M := M^{[r,\infty]}$ in terms
of $M^k := M^{[r,nk]}$, with $k\to \infty$. Then we follow the argument of Brown \cite{br87}.
\begin{itemize}  
\item First we note that $M$ is a directed partially ordered set and hence $| M|$  is contractible. 
\item $| M^{k+1} |$  is obtained from $| M^k |$  by adjoining cones over the subcomplexes $| M_{<f} |$ for each $f$ with
$h(f)=k+1$. By Lemma 3.7 and Remark 3.10 we know, that the subcomplexes $| M_{<f} |$  have the homotopy type of a
bouquet of $\big(h(S)- 1\big)$-spheres for k sufficiently large.
This shows that the embedding $| M^k | \subseteq | M^{k+1} |$  is homotopically trivial in all dimensions $< h(S).$
\item By Corollary 3.12 we know that the $| M^k |$ have cocompact skeleta.
\item The stabilizers, stab$_{G(S)}(f)$, of the vertices $f\in M$ -  in fact of all simplices - are of the form $G(S -
Sf)$. As rk$(S - Sf) < \text{rk}S$ the inductive hypothesis applies. The assumption that $M$ contains only injections $f$
with $h(f) \ge r$ implies now, that $fl \big({\rm stab}_{G(S)}(f)\big) \ge r - 1$ for each $f\in M$.
\end{itemize}

We can choose $r$ arbitrarily; if we choose $r \ge h(S) + 1$ the main results of \cite{br87} apply and it follows that $
fl\big(G_{\text{dia}}(S)\big) = h(S) - 1.$  This completes the inductive step.

In order to prove that $fl \big(G(S)\big) \ge fl\big(G_{\text{dia}}(S)\big)$ we note that $fl \big(G(S)\big) = fl\big(G_{\rm
tr}(S)\big)$, since $G_{\rm tr}(S)$ is of finite index in $G(S).$ Then we observe that $G_{\text{dia}}(S)$ is a normal subgroup of $G_{\rm
tr}(S)$ with $Q = G_{\rm tr}(S)/G_{\text{dia}}(S)$ finitely generated Abelian. As  $fl(Q) = \infty$ this implies
$fl\big(G_{\rm tr}(S)\big) \ge fl\big(G_{\text{dia}}(S)\big).$\hfill $\square$

\section{A lower bound for the finiteness length of pet$\mathbf{(S)}$ for a stack of orthants}
In this section we will show
 \begin{theorem} If $S$ is a stack  of orthants then $fl\big({\emph{pet}}(S)\big)\ge h(S) - 1.$
\end{theorem}
The steps to prove this lower bound of  $fl\big({\text{pet}}(S)\big)$ are similar to those in section 3 for
the corresponding pei-result. We will use a certain poset of injective pet-maps $f: S \to  S$ to form a simplicial
complex, and we will choose a diagonal subgroup of ${\text{pet}}(S)$ for the action on the complex.
However, the part concerning the finiteness length of the stabilizers of $f$ is more difficult here, because the set $(S - Sf)$ is generally not pet-isomorphic to a stack of orthants with lower rank
 (there are different parallelism classes of rank-$(n-1)$ germs in $S$ if rk$S = n)$.
So even if the stabilizers are isomorphic to pet$(S - Sf)$, there is no base for an induction argument.

In order to set up an inductive proof we need a version of Theorem 4.1, which makes the assertion not only
for stacks of orthants but also for stacks $S$ of paralell copies of a ``rank-n-skeleton'' of an orthant.  In combination with special
injective pet-maps $f$ (the ``super-diagonal'' maps), such a stack $S$ leads to a set $(S - Sf)$, which has the structure
of a stack of rank-$(n-1)$-skeletons.

\vskip 2mm 
\subsection{Stack of skeletons of an orthant.}~ Let $X$ be the canonical basis of the standard orthant $\mathbb N^N$.
Every orthant $L$ is of the form $a+\oplus_{y\in Y}\mathbb N y$, where $Y$ is a subset of $X$. $L$ carries the structure of a
simplex whose faces, indexed by the subsets $Z\subseteq Y$, are the suborthants $L_Z = a+\oplus_{z\in Z}\mathbb N z
\subseteq L.$ We refer to $L_Z$ as a \textit{rank-$k$ face} of $L$ if $| Z| =k.$ By the \textit{rank-$k$ skeleton} of $L$,
denoted $L^{(k)}$, we mean the union of all rank-$k$-faces of $L$. Thus the skeleta of $L$ form an ascending chain of
orthohedral set \[ \{ a\}  = L^{(0)} \subseteq  L^{(1)} \subseteq  \ldots \subseteq  L^{(k)} \subseteq  \ldots \subseteq  L^{({\text{rk}}L)}
 = L.\] Let $L^{(n)}$ be the rank-$n$ skeleton of a rank-$r$ orthant $L = a+\bigoplus_{ y\in Y}\mathbb N y$. Then
 $L^{(n)}$ is the union of $h(L^{(n)}) = \binom{r}{n}$ pairwise non-parallel rank-$n$ orthants.

Now we consider a stack $S$ of paralell copies of the rank-$n$-skeleton $L^{(n)}$ of an rank-$r$ orthant - in other words,
$S = R^{(n)}$ is the rank-n-skeleton of a stack $R$ of rank-$r$ orthants. We call each copy of $L^{(n)}$
in such a stack S a \textit{component of} $S$, and we write $c(S)$ for the number of components of $S$. Note that
$h(S)= c(S)\binom{r}{n}$.
The next proposition shows a lower bound for $fl\big(\text{pet}(S)\big)$, and the case $n=r$ yields the assertion of
Theorem 5.1.

\begin{proposition} If $S$ is a stack of rank-$n$ skeletons of an orthant then 
$fl\big(\emph{pet}(S)\big) \ge c(S) - 1$.
\end{proposition}

For later purpose in this section we consider the subset $\mathring{S}\subseteq S$ of all regular
points of $S$, which is defined as follows: If $S$ is an orthant, then $\mathring{S}$ is the image $t_S(S$) of $S$ under
the diagonal unit-translation; and if $S$ is a stack of rank-$n$ skeletons of an orthant, a point $p\in S$ is regular if $S$ contains a
maximal suborthant of rank equal to $S$, which contains $p$ as a regular point. The complement, denoted sing$(S) = S -
\mathring{S}$, is the set of all \textit{singular points} of $S$. 
 In the case when $S$ is a stack of orthants, we will also use the geometrically more suggestive notation $\partial S$ for 
 sing$(S)$. If $S = R^{(n)}$ is the $n$-skeleton of a stack of rank-$r$ orthants $L$, then sing$(S) = R^{(n-1)}$ and 
 $\mathring{S}$  has the canonical decomposition as the disjoint union of the regular points of the maximal orthants of $S$.
 By a \textit{component} of  $\mathring{S}$ we mean \mbox{$C \cap \mathring{S}$}, the intersection of $\mathring{S}$ with
 a component $C$ of $S$. Note that $c(S) = c(\mathring{S})$.
\begin{lemma} For the sets $S$ and $\mathring{S}$ the following holds
\begin{enumerate}\item $S$ and $\mathring{S}$ are pet-isomorphic. Hence  $\emph{pet}(S$) is isomorphic to
$\emph{pet}(\mathring{S})$.
\item $h\big(sing (\mathring{S})\big) = h\big(sing(S)\big)(r - n + 1)$, where $r$ is the rank of the stack $R$ with $S
= R^{(n)}$.
\end{enumerate}
\end{lemma}

\textit{Proof:} i) $S$ is the disjoint union of $\mathring{S}$ and $(S- \mathring{S})$. As each maximal orthant of $(S-
\mathring{S})$ is parallel to a subortant of $\mathring{S}$, the assertion follows from the pet-normal form.
\\
ii) Since $\mathring{S} = R^{(n)} - R^{(n-1)} $, sing$(\mathring{S})$ is the disjoint union of $h(R^{(n)}) \cdot n$
rank-$(n-1)$ orthants. So $h(\text{sing}(\mathring{S})) = h(R^{(n)})n$. For the height of $S$ and sing$(S)$ we have
$h(S)=h(R^{(n)}) = c(S)\binom{r}{n}$ and $h(\text{sing}(S))=h(R^{(n-1)})= c(S)\binom{r}{n-1}$. As $\binom{r}{n} n =
\binom{r}{n-1} (r-n+1)$, we get $h(R^{(n)}) n = h(R^{(n-1)})(r-n+1)$.
\hfill$\square$
\vskip 2mm 
\subsection{Reduction to the diagonal subgroup.}~ From now on we assume that $S$ is a stack of rank-$n$ skeletons of an
orthant. Since $S$ and $\mathring{S}$ are pet-isomorphic, it suffices to show Proposition 4.2 for the set
$\mathring{S}$, which is more suitable for some parts of the proof.
As noted above $\mathring{S}$ is canonically in pet-normal form. 
In particular, every maximal germ of S (or $\mathring{S}$) is represented by a unique
maximal orthant of $\mathring{S}$. Thus we can conceptually simplify matters by replacing the set of all maximal germs,
max$\Gamma^\ast(S)$ = max$\Gamma^\ast(\mathring{S})$, by the set of the canonical representatives
max$\Omega^\ast(\mathring{S})$, the set of all maximal orthants of $\mathring{S}$.
\\
Let $M_{\text{dia}} (\mathring{S})$ denote the monoid of all diagonal
pei-injections of $\mathring{S}$ introduced in Section 3.2. $M_{\text{dia}} (\mathring{S})$ is a submonoid of
$M_{\text{tr}} (\mathring{S})$, the translation submonoid of $M(\mathring{S})$.
Its elements f have the property that they induce, for each
$L \in \text{max}\Omega^\ast(\mathring{S})$, a diagonal translation $\tau_{(f,L)}: \langle L \rangle \to \langle L
\rangle$. 
 Now we consider the submonoid $M^{\text{pet}}_{\text{sdia}}(\mathring{S})\subseteq M_{\text{dia}} (\mathring{S})$
 consisting of all diagonal pet-injections $f: \mathring{S}\to \mathring{S}$ which satisfy the additional \textit{super-diagonality} condition:

\begin{tabular}{@{}lP{14.2cm}@{}}
(4.1) & \textit{When two maximal orthants $L, L' of  \mathring{S}$ are contained in the same component of  $S$,
then the diagonal translations $\tau_{ (f,L)}$ and $\tau_{ (f,L')}$  have the same translation length.}
\end{tabular}
 
The restriction of the homomorphism (3.3) of Section 3.2 to $M^{\text{pet}}_{\text{sdia}}(\mathring{S})$ can thus be interpreted as a
map

\begin{tabular}{@{}lP{14.2cm}@{}}
(4.2) &	$\lambda : M^{\text{pet}}_{\text{sdia}}(\mathring{S})\twoheadrightarrow \bigoplus_{ C\in Comp(\mathring{S})}\mathbb Z = \mathbb
Z^{c(S)}$,
\end{tabular}

which associates to each super-diagonal pet-injection $f$ the translation length $\lambda 
(f,C)$ on each component $C$ of $\mathring{S}$.

The group of all invertible elements of $M^{\text{pet}}_{\text{sdia}} (\mathring{S})$ is the \textit{super-diagonal pet-group} 
pet$_{\text{sdia}} (\mathring{S})$. 

Let be pet$_{\text{tr}} (\mathring{S})$ the group of all invertible elements of $M_{\text{tr}} (\mathring{S})$.
It is a subgroup of $\text{pet}(\mathring{S})$, which has finite index in $\text{pet}(\mathring{S})$.
Analogous to (3.2) in Section 3.2 there is a homomorphism

\begin{tabular}{@{}lP{14.2cm}@{}}
(4.3)  & 	$\kappa : \text{pet}_{\text{tr}}(\mathring{S})\twoheadrightarrow \bigoplus_{ L\in \max
\Omega^\ast(\mathring{S})} Tran(\langle L \rangle ),$ given by $\kappa (g) = \bigoplus_{ L\in \max
\Omega^\ast(\mathring{S})}\tau (g,L)$
\end{tabular}

which associates to each pet-injection $g\in \text{pet}_{\text{tr}}(\mathring{S})$ the translation length $\tau 
(g,L)$ on each maximal orthant $L$ of max$\Omega^\ast(\mathring{S})$. We observe that a permutation $g\in
\text{pet}_{\text{tr}}(\mathring{S})$ is in $\text{pet}_{\text{sdia}} (\mathring{S})$, if and only if the translations 
$\tau_{ (g,L)}$ are diagonal for each $L$ and its translation length constant as $L$ runs through the maximal orthants of
a component $C$ of $\mathring{S}.$ 
\\
Given a component $C$ of $\mathring{S}$, we consider the set $\Lambda(C):=\text{max}\Omega^\ast({C})$ of all $h(C) =
\binom{r}{n}$ rank-$n$ orthants of $C$. For each orthant $L\in \Lambda(C)$, we
write $Y(L)$ for its canonical basis. The translation $\tau_{ (g,L)} :\langle L \rangle \to\langle L\rangle$ has the
canonical decomposition into the direct sum of translations $\tau^y_{(g,L)}$ in the directions $y\in Y(L)$, and we write
$l^y(g, L)\in \mathbb Z$ for the corresponding translation lengths. 
\\
Hence, for $g\in \text{pet}(\mathring{S})$ to be super-diagonal means, that the numbers  $l^y(g, L)\in\mathbb Z$  coincide 
for all pairs in $P(C) := \{ (y, L)\mid y\in L\in \Lambda(C) \}$ 
- and this is so for all components $C$. Hence, associating to $g$ the sequence \mbox{$\big( l^y(g, L) - l
^{y'}(g,L')\big)_{( i(C),C)}$}~, with $i(C)$ running through all pairs $\big((y, L), (y', L')\big) \in  P(C)$, and
$C$ through the components of $\mathring{S}$, exhibits the super-diagonal pet-group pet$_{\text{sdia}} (\mathring{S})$
as the kernel of a homomorphism of $\text{pet}_{\text{tr}}(\mathring{S})$ into
a finitely generated Abelian group. It is well known that in this situation 
$fl\big(\text{pet}_{\text{tr}}(\mathring{S})\big) \ge fl\big(\text{pet}_{\text{sdia}} (\mathring{S})\big)$. Since
$\text{pet}_{\text{tr}}(\mathring{S})$ has finite index in $\text{pet}(\mathring{S})$, we have
$fl\big(\text{pet}(\mathring{S})\big) = fl\big(\text{pet}_{\text{tr}}(\mathring{S})\big)$, hence
$$fl\big(\text{pet}(\mathring{S})\big) \ge fl\big({\text{pet}}_{\text{sdia}} (\mathring{S})\big).$$ 
The proof of Proposition 4.2 is thus reduced to a proof of $fl\big({\text{pet}}_{\text{sdia}} (\mathring{S})\big) =
c(S) - 1.$ To show this, we follow the arguments in the proof of the corresponding pei-result: 
$fl\big({\text{pei}}_{\text{dia}}(S)\big) = h(S) - 1$, where $S$ was a stack of $h(S)$ orthants of rank $n$. 
In the present situation, where $\mathring{S}$ is the set of regular points of the $n$-skeleton of the stack $R$ of $h(R)$ 
orthants, the components $C$ of $\mathring{S}$ have to take over the role previously played by the orthants $L$ of the stack $S$.
Correspondingly we now have to work with the multiplicative submonoid $\text{mon} (T)\subseteq M^{\text{pet}}_{\text{sdia}} (\mathring{S})$ 
freely generated by the set $T$ of all super-diagonal unit-translations $t_C: C\to C$ as $C$ runs through the
components of $\mathring{S}$, where each $t_C=\prod_{L\in \Lambda(C)}t_L$ is the composition of the diagonal unit-translations
$t_L$ defined in Section 3.2. As at the end of Section 3.2 we use the action of $\text{mon} (T)$ by left multiplication to
endow $M^{\text{pet}}_{\text{sdia}} (\mathring{S})$ with a partial ordering; and we observe that this partial ordering is directed.
\vskip 2mm 
\subsection{Maximal elements $\mathbf{< f}$ in $\mathbf{M^{\text{pet}}_{\text{sdia}} (\mathring{S})}$.}~  
To adapt notation to the one used in the corresponding pei-situation in Section 3, we write $\Lambda$  for the set of all
components $C$ of $\mathring{S}$, and $\partial C := C - Ct_C$  for each component $C\in \Lambda$ .
Note that $C$ is the disjoint union of $h(C) = \binom{r}{n}$ rank-$n$ orthants
-- using the notation of Section 4.1 one for each $n$-element set $Z\subseteq Y$.
Hence $h(\partial C) = n \binom{r}{n}.$ We are still in the situation that all maximal orthants of $\mathring{S}$ have the same
 finite rank $n = rk S.$  And given $f\in M^{\text{pet}}_{\text{sdia}} (\mathring{S})$ we write
\[	M_{ <f}   = \{ a\in M^{\text{pet}}_{\text{sdia}} (\mathring{S})\mid a <f  \} , ~ 	M_{\le f} = \{ a\in M^{\text{pet}}_{\text{sdia}} (\mathring{S})\mid a\le f\} ,\]  
for the \textit{``open resp. closed cones below $f$ ''}, aiming to understand the set of all maximal
elements of $M_{<f}$.

\begin{lemma} Let $b$ be a maximal element of  $M_{<f}$. Then there is a unique component $C$ of $\mathring{S}$ with the 
property that $f = t_Cb$, and $h(f) = h(b) + n$. Furthermore, $b$ is given as the union $b = b'\cup b''$, where b': 
$\partial C\to (\mathring{S} - \mathring{S}f)$ is a pet-injection, and $b'': (\mathring{S} - \partial C) \to 
\mathring{S}f$ is the restriction $t_C^{-1}f|_{(\mathring{S} - \partial C)}$.
Converseley, if $c': \partial C\to (\mathring{S} - \mathring{S}f)$ is an arbitrary pet-injection distinct to
$b'$, then the union $c = c'\cup b''$ is a maximal element of $M_{< f}$ distinct to $b$.
\end{lemma}

\textit{Proof.} See argument in Lemma 3.4. \hfill$\square$

\begin{lemma} Let $B\subseteq M_{<f}$  be a finite set of maximal elements of $M_{<f}$  . Then the following conditions are equivalent:
\begin{enumerate} \item The elements of $B$ have a common lower bound $\delta $ in $ M_{<f}.$
\item For every pair $(b,b')\in B\times B$, with $b\not= b'$ and $tb =f =t'b'$ for super-diagonal unit-translations $t$,
$t'$, we have
\begin{enumerate}\item		 $t\not= t'$,  and 
\item		  $b(\partial C )\cap b'(\partial C') = \emptyset$,
 where $C$ resp. $C'$ are the components of $\mathring{S}$ on which $t$ resp. $t'$ acts non-trivially.
\end{enumerate}
\end{enumerate}
\end{lemma}
\textit{Proof.} For each $b\in B$ we have a super-diagonal unit-translation $t_b\in T$ with $t_b b = f$, and we put

\begin{tabular}{@{}lP{14.2cm}@{}}
(4.4)& $t_B := \prod_{ b\in B}t_b$.
\end{tabular}

By assumption (a) the components $C_b$, on which $t_b $ acts non-trivially, are pairwise disjoint. 
Thus $|B|$ $\le c(S)$, and $S$ decomposes in the disjoint union  $S  =
\left(\bigcup_{b\in B}C_b \right) \cup  S'$.
We define the pet-injection $\delta_B: S\to S$ as follows:
$$ \delta_B := \begin{cases}
{t_b}^{-1} f & \text{on each } C_b t_b \\ 
b & \text{on the complements } \partial C_b = C_b - C_b t_b \\ 
f & \text{on } S'. 
\end{cases} $$
To show that $\delta_B$ is a common lower bound, see arguments in Lemma 3.5.  \hfill$\square$

\begin{lemma} In the situation of Lemma 4.5 we have for the lower bound
$\delta_B$ defined in the proof:
\begin{enumerate}\item $\delta_B$ is, in fact, a largest
common lower bound of the elements of $B$
\item $h(\delta_B ) \ge h(f) - h(S)n$.
\end{enumerate}
\end{lemma} 	
\textit{Proof.} For (i) see argument in Lemma 3.6. For (ii) we use the translation  $t=
\Pi_{ b\in B}t_b$ of (4.4) which satisfies $t\delta_B  = f$  and reads:
\[ 
\begin{array}{lclr}
h(\delta_B ) = h(f) - h(t ) &=& h(f) - |B| \cdot h(t_C) &\\
&=&  h(f) - |B| \cdot h(C)n &\\
&\ge& h(f) - c(S)h(C)n = h(f) - h(S)n.&
\end{array}
\]\hfill$\square$

\subsection{The simplicial complex of $\mathbf{M^{\text{pet}}_{\text{sdia}} (\mathring{S})}$.}~ We consider the simplicial complex
$|M^{\text{pet}}_{\text{sdia}} (\mathring{S})|$, whose vertices are the elements of $M^{\text{pet}}_{\text{sdia}} (\mathring{S})$ and whose chains of
length $k,  a_0  < a_1  < \ldots < a_k$, are the $k$-simplices. As the partial ordering on $M^{\text{pet}}_{\text{sdia}} (\mathring{S})$
is directed, $|M^{\text{pet}}_{\text{sdia}} (\mathring{S})|$  is contractible.

In this section we aim to prove    
\begin{lemma} If $h(f)\ge2\cdot \emph{rk}S\cdot h(S)$ then $|M_{<f}|$  has the homotopy type of a bouquet of
$\big(c(S) - 1\big)$-spheres.
\end{lemma}
 The first step towards proving Lemma 4.7 is to consider the covering of  $|M_{<f}|$  by the subcomplexes $|M_{\le b}|$,
 where $b$ runs through the maximal elements of  $M_{<f}$. We write $N(f)$ for the nerve of this covering. Lemma 4.6(i) asserts that all finite intersections of such subcomplexes $|M_{\le b}|$  are again cones and hence contractible. It is a well known fact that in this situation the space is homotopy equivalent to the nerve of the covering. Hence we have \[  |M_{<f}| ~\text{is homotopy equivalent to the nerve}~ N(f),\]
and it remains to compute the homotopy type of~ $N(f).$

The next step - replacing the nerve $N(f)$ by the combinatorial complex $\Sigma(f)$ - follows the arguments in Section 3:
We find that the set of vertices of $\Sigma(f)$ is the disjoint union  $A = \bigcup_{ C\in \Lambda} A_C $, where $A_C $
stands for the set of all pet-injections  $a: \partial C\to (\mathring{S} - \mathring{S}f )$; 
and the $p$-simplices of $\Sigma(f)$ are the sequences $(a_C)_{C\in \Lambda'}$, where $\Lambda ' $ is a $p$-element subset
of $\Lambda$ whose entries $a_C\in A_C $ satisfy the condition

 (4.5) \quad \textit{the intersections of the images}~ $ a_C(\partial C), C\in \Lambda '$, ~\textit{are pairwise
 disjoint.}

The homotopy type of $\Sigma(f)$ can again be computed by Lemma 3.8., which we apply to the $1$-skeleton  $\Gamma (f)$ of 
$\Sigma(f)$, viewed as a $c(S)$-colored graph $\Gamma (f)_{c(S)}$.
At the end it remains to prove 
\begin{lemma} If  $h(f)
\ge 2\cdot\emph{rk}S\cdot h({S})$  then $\Gamma (f)_{c(S)}$ satisfies the assumptions of Lemma 3.8.
\end{lemma}
\textit{Proof.} Let $n := \text{rk}S$ and $h := c(S)$. Assumption (1) is a consequence of assumption (2) except in the trivial case $h = 1.$ 

To prove (2) we fix $C\in \Lambda$  and consider a set of $2(h - 1)$ elements $F\subseteq A - A_C$ . We have to show that
there are two elements  $a, b\in A_C$  with the property that for each $d\in F, d: \partial C_d \to  (\mathring{S} -
\mathring{S}f), im(d) = d(\partial C_d)$ is disjoint to both $a(\partial C)$ and $b(\partial C).$ In other words: there
are two pet-injections

(4.6)	\quad	$a, b:  \partial C  \to  (\mathring{S} - \mathring{S}f) - \big(\bigcup_{ d\in F} im(d)\big).$
 
For this it suffices to compare the height function - i.e., the number of rank-$(n-1)$ germs - of domain and target. 
Clearly, $h\big(a(\partial C)\big) = h(\partial C) = nh(C)$, and the same applies to every vertex of $A.$ Hence 
$h\big(\bigcup_{ d\in F} im(d)\big) \le 2(h - 1)nh(C)$. By assumption $h(\mathring{S} - \mathring{S}f) \ge
2nh({S})$, and so the target orthohedral set has height at least $2nh({S}) - 2(h - 1)nh(C) = 2n h(C)$, which is more than at 
least twice the height $h(\partial C) = nh(C)$ of the domain when $h(C)$ is positive. Moreover, by Lemma 4.10 below, the
set $(\mathring{S} - \mathring{S}f)$ is pet-isomorphic to a stack of copies of $\partial C$. In this situation one observes that the two different pet-injections required in (4.6) certainly do exist. 
This proves the Lemma 4.8 and hence Lemma 4.7. \hfill$\square$

\begin{remark}
\emph{
By the same argument as in Remark 3.10, the assertion of Lemma 4.7 holds true if $M_{<f}$ is replaced with the subset
$M_{r,f} :=\{a \in M^{\text{pet}}_{\text{sdia}} (\mathring{S}) \mid h(a) \geq r ~ and ~ a < f \}$ and $f$ satisfies the
additional condition $h(f) \geq r + h(S)$.}
\end{remark}
\vskip 2mm 
\subsection{Stabilizers and cocompact skeletons of $\mathbf{|M^{\text{pet}}_{\text{sdia}} (\mathring{S})|}$.}~ Here we consider the monoid
$M^{\text{pet}}(S)$ of all pet-injections endowed with the height function $h: M^{\text{pet}}(S)\to \mathbb Z$ inherited from $M(S)$ and the
pet$(S)$-action induced by right multiplication. Our main interest, however, is the
super-diagonal submonoid $M^{\text{pet}}_{\text{sdia}} (\mathring{S})\subseteq M^{\text{pet}}(S)$ acted on by the super-diagonal pet-group
${\text{pet}}_{\text{sdia}} (\mathring{S}).$ 
\begin{lemma} \begin{enumerate}\item $f,f' \in M^{\emph{pet}}(S)$ are in the same
\emph{pet}$(S)$-orbit if and only if $(S - Sf )$ and $(S - Sf ')$ are 	pet-isomorphic.
\item Let $S$ be a stack of copies of $L^{(n)}$, where $L$ is a rank-$r$ orthant. If $f\in M^{\emph{pet}}_{\emph{sdia}}
(\mathring{S})$ with $h(f) > 0$, then $h(f)$ is a multiple of $(r-n+1)\binom{r}{n-1}$ and $(\mathring{S} - \mathring{S}f)$
is pet-isomorphic to a stack of $h(f) / \binom{r}{n-1}$ copies of  $L^{(n-1)}$.
\item Two elements $f,f '\in M^{\emph{pet}}_{\emph{sdia}} (\mathring{S})$ with $h(f) = h(f ') > 0$ are in the same
$\emph{pet}_{\emph{sdia}} (\mathring{S})$-orbit.
\end{enumerate}
\end{lemma}
\textit{Proof.} The proof of i) is analogous to the pei-version in Section 3.\\  
ii) The key here is a pet-version of Lemma 3.2i). The set of germs $\Gamma^{ n-1}~(\mathring{S})$ decomposes into its
parallelism classes, and as these are pet$(\mathring{S})$-invariant, the height function 
$h: M^{\text{pet}}_{\text{sdia}} (\mathring{S})\to \mathbb Z$ can be written as the sum of functions 
$h_Y: M^{\text{pet}}_{\text{sdia}} (\mathring{S})\to \mathbb Z$, with $Y$ running through all $(n-1)$-element subsets of the canonical basis 
of $L$, that count the number of germs in $\Gamma^{ n-1}(\mathring{S} - \mathring{S}f)$ parallel to $\langle Y\rangle .$

We need the $h_Y $-version of Lemma 3.2i), asserting that we have for all rank-$(n-1)$ faces $Y$ of $L$
  
\begin{tabular}{@{}lP{14.2cm}@{}}
(4.7) & $h_Y (f)  =  h_Y \big(A\cap (Af)^c\big) - h_Y (A^c \cap Af)$, \textit{for every orthohedral subset} $A\subseteq S$ \\
      & \textit{with} $\text{rk}A^c  < n$.
 \end{tabular}
  
The proof is a straightforward adaptation of the one in Section 3.1 and can be left to the reader.

We can refine (4.7) by exhibiting $A$ as the disjoint union of rank-$n$ orthants $K_i$  on which $f$ acts by
(super-diagonal) translations. Since $f$ is super-diagonal, the corresponding translation lenghts $\lambda_{C(i)}$ depend
only on the component $C(i)$ of $\mathring{S}$ containing $K_i$. One observes that $f(K_i )$ is
contained in the uniquely defined maximal orthant of $\mathring{S}$ containing $K_i$. This has the consequence that for
$i \not=  j,  K_i  \cap (K_jf)^c  = K_i$   and $K_i^c\cap K_jf = K_jf,$ from which one finds

\begin{tabular}{@{}lP{1.2cm}@{}cl}
(4.8)  &   $h_Y (f)$ &=& $h_Y  \big(\bigcup_i (K_i \cap (K_i f)^c\big)- h_Y  \big(\bigcup_i ( K_i^c \cap K_i f)\big)$\\
       &             &=& $\sum_i h_Y(K_i \cap (K_i f)^c -  h_Y( K_i^c \cap K_i f) = \sum_i \lambda_{C(i)}$.
\end{tabular}

Clearly, for each
component $C$ of $S$, $sing(C)$ contains exactly one orthant parallel to $\langle Y\rangle$. By Lemma 4.3ii) this orthant
is parallel to a face of exactly $(r - n + 1)$ maximal orthants $K_i$  in $C$, and each of them gives rise to a summand
$\lambda_{C(i)}$ .
Hence summation over all $K_i$  contained in a single component $C$ of $S$ yields  $\lambda_C(r - n + 1)$. And summation over all I, finally,

 $h_Y (f)  =  \lambda (r - n + 1)$ , where $\lambda$  is the sum of $\lambda_C$,  with $C$ running through all components of $S$.

This shows, in particular, that $h_Y (f)$ is intependent of $Y$. As we are assuming that $h(f) > 0$ it follows that 
$\lambda  > 0$  and  $h(f) = \lambda \cdot (r - n + 1)\cdot \binom{r}{n-1}$.

It follows that $(\mathring{S} - \mathring{S}f)$  is pet-isomorphic to disjoint union $S'\cup S''$, where $S'$ is a stack of 
$h_Y (f) = \lambda (r - n + 1)$ copies of $L^{(n-1)}$ and $S''$ a subset of rank $< n - 1$. As $\lambda  > 0$, $S'$
contains at least one copy of $L^{(n-1)}$. In this situation $S'$ contains orthants parallel to any given maximal orthant of
$S''$. In view of the the pet-normal form of  $S'\cup S''$  it follows that  $S'\cup  S''$ is pet-isomorphic to $S'$, i.e., 
to a stack of    $\lambda (r - n + 1)$ copies of $L^{(n-1)}$.

iii) Part ii) shows that $h(f) = h(f ') > 0$ implies that $\mathring{S} - \mathring{S}f$ and $\mathring{S} - \mathring{S}f
'$ are pet-isomorphic. Hence, by assertion i), there is a pet-permutation $g\in {\text{pet}}(\mathring{S})$ with $f ' =
fg$. The assumption that $f$ and $f '$ are in $M^{\text{pet}}_{\text{sdia}} (\mathring{S})$ implies that $g\in
{\text{pet}}_{\text{sdia}} (\mathring{S}).$\hfill$\square$

 \begin{corollary} \begin{enumerate} \item The stabilizer of  $f\in M^{\emph{pet}}_{\emph{sdia}} (\mathring{S})$ in
 ${\emph {pet}}_{\text{sdia}} (\mathring{S})$ is isomorphic to ${\emph{pet}}(\mathring{S} - \mathring{S}f)$
\item For every number $r, s\in \mathbb N_0$, the simplicial complex of $M^{[r,s]}:= \{ f\in
M^{\emph{pet}}_{\emph{sdia}} (\mathring{S})\mid r\le h(f)\le s\}$  is cocompact under the \emph{pet}$_{\emph{sdia}}
(\mathring{S})$-action.
\end{enumerate}
\end{corollary}

\textit{Proof.} i) Same reasoning as the proof of Corollary 3.12. 
\\ii) Lemma 4.10iii) asserts that ${\text{pet}}_{\text{sdia}} (S)$ acts
cocompactly on the vertices of a given height in the  simplicial complex $| M^{\text{pet}}_{\text{sdia}} (S)|$. Just as in the proof 
of Corollary 3.12ii) this yields the claimed assertion.
\hfill$\square$
\vskip 2mm 
\subsection{The conclusion.}~ It is an elementary observation that right action of $g\in {\text{pet}}(S)$ on
$M^{\text{pet}}(S)$ fixes an element $f\in M^{\text{pet}}(S)$ if and only if $g$ restricted to $f(S)$ is the identity. In other words, the stabilizer of a vertex
$f\in M^{\text{pet}}(S)$ is isomorphic to \mbox{pet$(S - Sf)$}. This will be crucial for the inductive step in the
following inductive

\textit{Proof of Proposition 4.2}. In Section 4.2 was already proved $fl\big(\text{pet}(\mathring{S})\big) \geq
fl\big(\text{pet}_{\text{sdia}}(\mathring{S})\big)$, so it suffices to show
$fl\big({\text{pet}}_{\text{sdia}}(\mathring{S})\big) \ge c(S) - 1$. We will argue by induction on $n = \text{rk}S$.
If $n = 1$, then $h(S)=c(S)\cdot r$, and the group pet$_{\text{sdia}}(\mathring{S})$ is the Houghton group on $h(S)$ rays.
By Brown \cite{br87} this implies that $fl\big({\text{pet}_{\text{sdia}}}(\mathring{S})\big) \ge h(S) - 1 \ge c(S) - 1$.
This establishes the case $n = 1$ of the induction.
      
Now we assume $n > 1$. By induction we can assume that $fl\big({\text{pet}}(S')\big) \ge c(S') - 1$ holds for
every stack S' of copies of a rank-$(n-1)$ skeleton of an orthant $L$. To prove the inductive step we start with restricting
attention to the subgroup pet$_{\text{sdia}} (\mathring{S})$ acting on the super-diagonal monoid
$M^{[u,v]} = \{ f\in  M^{\text{pet}}_{\text{sdia}} (\mathring{S})\mid u\le h(f)\le v\} , u,v\in  \mathbb N$. By
Lemma 4.10ii) the $h(f)$  is a multiple of $s := (r - n + 1)\binom{r}{n-1}$. So we fix a lower bound $u = sk_0, k_0\in 
\mathbb N$, and consider the filtration of  $M:= M^{[u, \infty]}$ in terms of $M^k = : M^{[u,sk]},$ with $k\to \infty.$
Then we argue as follows
\begin{itemize}
\item First we note that $M$ is a directed partially ordered set and hence $|M|$  is contractible. 
\item $| M^{k+1} |$  is obtained from $| M^k |$  by adjoining cones over the subcomplexes $| M_{<f} |$ for each $f$
with $h(f)=k+1$. By Lemma 4.7 and Remark 4.9 we know, that the subcomplexes $| M_{<f} |$  have the homotopy type of a
bouquet of $\big(c(S)- 1\big)$-spheres for k sufficiently large. This shows that the embedding $|M_k| \subseteq |M_{k+1}|$
is homotopically trivial in all dimensions $< c(S)$.
 \end{itemize}
\begin{itemize}
\item By Corollary 4.11 we know that each $| M^k|$  has cocompact skeleta.
\item The stabilizer of the $f\in M$ under the action of pet$_{\text{sdia}} (\mathring{S})$ on $M$ coincides with
pet$(\mathring{S} - \mathring{S}f)$ by Corollary 4.11i). Lemma 4.10ii) asserts that if $h(f) > 0$ then $(\mathring{S} -
\mathring{S}f)$ is pet-isomorphic to a stack of copies of the rank-$(n-1)$ skeleton of an orthant.  The stack height here is
 $c(\mathring{S} - \mathring{S}f) = h(f)/\binom{r}{n-1}$. We can choose u arbitrarily; if we choose $u = \big(c(S) +  1)
 \binom{r}{n-1}$ the inductive hypothesis together with the assumption that $h(f)\ge u$ yields
 \[fl\big({\text{pet}}(\mathring{S} - \mathring{S}f)\big) \ge c(\mathring{S} - \mathring{S}f) - 1 = h(f) \left/
 \binom{r}{n-1} - 1 \right. \ge c(S), ~\text{for all}~ f\in M.\]
\item The main results of \cite{br87} now establishs $fl\big(\text{pet}_{\text{sdia}}(\mathring{S})\big) \ge c(S) - 1.$
This completes the inductive step.
\hfill$\square$
\end{itemize}

\section{The upper bounds of $\mathbf{fl \big({\text{pet}}(S)\big) }$ when $\mathbf{S \subseteq  \mathbb N^N}$}
\subsection{More structure at infinity.}~ Here we assume, for simplicity, that our orthohedral sets $S$ are contained in
$\mathbb N^N$. By Corollary 1.5 this is not a restriction for the pei-group pei$(S)$, and it is a basic
special case for the  pet-group pet$(S)$:

Given an element $x\in X$ (i.e. a coordinate axis), we write $\Gamma_{x}^1(S)$ for the set of all germs of rank-$1$
orthants of $S$ parallel to $\mathbb N x$. We have a canonical embedding $\kappa : \Gamma_{x}^1(S)\to  \mathbb N^{ N-1}$
defined as follows: Each $\gamma \in \Gamma_{x}^1(S)$ is represented by a unique maximal orthant $L\in \Omega^1(S)$; we
delete the $x$-coordinate of the base point of $L$ and put $\kappa (\gamma )$ to the remaining coordinate vector. We write
$\partial_xS$ for the image $\kappa \big(\Gamma_{x}^1(S)\big)$, and we will often identify $\Gamma_{x}^1(S)$ with
$\partial_xS$ via $\kappa .~ \partial_xS$ can be viewed as the boundary of $S$ \textit{at infinity in direction}~ $x$.
\begin{lemma}For $\partial_xS$ the following hold.
\begin{enumerate} 
\item $\partial_xS$ is an orthohedral subset of $\mathbb N^{N-1}.$
\item For very rank-$(k-1)$ orthant $L\subseteq \partial_xS$ there is a unique rank-$k$ orthant $L'\subseteq S$, which
is maximal with respect to the property that for each point of $p\in L~ \kappa^{ -1}(p)$ is  represented by a suborthant of $L.$
\end{enumerate}
\end{lemma}

\textit{Proof.} Easy. \hfill$\square$

\subsection{Short exact sequences of pet-groups.}~ From now on we assume that $S = \bigcup_j S_j \subseteq  \mathbb N^m$
where $m$ is minimal and $S$ is in pet-normal form as defined after Proposition 1.5. Given $x\in X$ arbitrary we note that
$S$ is the disjoint union $S = S(x)\cup S^{\bot}(x)$, where $S(x)$ collects the stacks $S_j$ which contain a rank-$1$ orthant parallel to $\mathbb N
x$, and $S^{\bot}(x)$ the stacks $S_j$ which are perpendicular to $x$. We note that $\partial_x S = \partial_x S(x)$, and
we have an obvious projection $\pi_x: S(x) \twoheadrightarrow \partial_x S$. Moreover, there is a canonical injection 
$\sigma_x: \partial_x S \to  S(x)$ which maps each germ $\gamma \in \partial_x S$ to the base point of the unique maximal 
rank-$1$ orthant representing $\gamma$, and is right-inverse to  $\pi_x: S(x) \twoheadrightarrow \partial_x S.$

As the action of pet$(S)$ on the $\Omega^1(S)$ preserves directions it induces, for each coordinate axis $x\in X$, an
action on $\Gamma_{x}^1 (S) = \partial_x S$, and one observes that this is an action by pet-permutations. This yields an
induced homomorphism \mbox{$\vartheta_x : \text{pet}(S)\to \text{pet}(\partial_x S)$}. The kernel of $\vartheta_x$  is the set of
all pet-permutations fixing all rank-$1$ germs parallel to $x$. And we note the following:

\begin{enumerate}
\item
\mbox{$\sigma_x: \partial_x S \to  S(x)$} induces an embedding of pet$(\partial_x S)$ as a subgroup of
 pet$\big(S(x)\big)$, which splits the surjective homomorphism \mbox{$\vartheta_x : {\text{pet}}(S(x)) \twoheadrightarrow 
 {\text{pet}}(\partial_x S)$} induced by $\pi_x$.
\item  Every pet-permutation $g\in {\text{pet}}\big(S(x)\big)$ extends to a pet-permutation of $S$ by the identity on $S^{\bot}(x)$. 
This exhibits ${\text{pet}}\big(S(x)\big)$ as a canonical subgroup of pet$(S)$. Even though we do not have pet$(S)$ acting on
$S(x)$, we do have that the surjective homomorphism \mbox{$\vartheta_x : {\text{pet}}(S)\twoheadrightarrow \text{pet}(\partial_x
S)$} splits by the embedding ${\text{pet}}(\partial_x S)\le {\text{pet}}\big(S(x)\big)\le {\text{pet}}(S).$
\end{enumerate}
Summarizing we have
\begin {proposition}
${\emph{pet}}(\partial_x S)$ is a retract both of ${\emph{pet}}(S)$ and of ${\emph{pet}}\big(S(x)\big)$. In other words,
we have split exact sequences  $1 \to  K \to  {\emph{pet}}(S) \to {\emph{pet}}(\partial_x S) \to  1$  and  $1 \to  K^\dag
\to  {\emph{pet}}\big(S(x)\big) \to {\emph{pet}}(\partial_x S) \to  1.$ \hfill $\square$
\end{proposition}

\subsection{An upper bound of the finiteness length of $\mathbf{pet(S)}$.}~  To deduce an upper bound for the finiteness
lengths of the pet-groups we need the following elementary lemma which was overlooked in [BE75]; the first occurrence in
the  literature we are aware of is [Bu04].
\begin{lemma}
Let $G$ be a group. If a subgroup $H\le G$ is a retract of $G$ then  $fl(H) \ge fl(G).$
\end{lemma}
\textit{Proof.} The assertion  $fl(G)\ge s$ is equivalent to saying that on the category of $G$-modules the homology functors  
$H_k(G; - )$ commute with direct products for all $k<s.$ That this is inherited by retracts follows from the fact that 
$H_k( - ; - )$ is a functor on the appropriate category of pairs $(G,A)$, with $G$ a group and $A$ a $G$-module. 
\hfill$\square$

If $S \subseteq \mathbb N^N$ is an orthohedral set of rank $\text{rk}S = n$ in pet-normal form, then
$\Omega_0^\ast(\mathbb N^m)$ is canonically bijective to the set $P(X)$ of all subsets of $X$. Hence we can view the
height function (1.1) of section 1.4 as map \mbox{$h_S:
P(X) \to \mathbb N_0$}, and organize stacks of maximal orthants of $S$ as follows: For every subset $Y\subseteq X$ we have
the (possibly empty) stack $S(Y) \subseteq S$ of $h_S(Y)$ orthants parallel to the orthant $\langle Y \rangle$ defined by
$Y$.
\\For each $(n-1)$-element subset $Y\subseteq X$ we consider the link ${\text{Lk}}(Y)$ of $Y$ in $S_\tau$, by this we mean the set
of all $n$-element sets $Y'\subseteq Y$ with the property, that $\langle Y' \rangle\subseteq S_\tau$, noting that $\langle
Y' \rangle \in {\rm max} \Omega_0^\ast(S_\tau)$. Then we put $S({\text{Lk}}(Y)) \subseteq S$ to be the union of the stacks $S(Y')$
with $Y'$ running through ${\text{Lk}}(Y)$. The height, $h\big(S\big({\text{Lk}}(Y)\big)\big)$, is the sum of all
stack-heights $h_S(\langle Y' \rangle)$ as $Y'$ runs through ${\text{Lk}}(Y)$.

\begin{theorem}
If $S \subseteq  \mathbb N^m$ is orthohedral of rank $n$ in pet-normal form, then each indicator $(n-1)$-subset $Y
\subseteq X$ with non-empty link ${\emph{Lk}}(Y)$ imposes an upper bound $fl\big(\emph{pet}(S)\big)\leq
h\big(S\big({\emph{Lk}}(Y)\big)\big)-1$.
\end{theorem}

\textit{Proof.} We choose any $y\in Y$ and consider the projection \mbox{$\pi_y: S \twoheadrightarrow
\partial_y S\cup\{\emptyset\}$}, where the symbol $\{\emptyset\}$ is the image of $S-S(y)$. Proposition 5.2 asserts, that
pet$(\partial_y S)$ is a retract of pet$(S)$. We have $\text{rk}\partial_yS=\text{rk}S-1$; in fact, $\partial_yS$ is the
disjoint union of stacks $S(Z)$, with $Z$ running through all subsets of $X$ avoiding y and satisfying $\langle Z\cup \{y\}\rangle\in
{\rm max} \Omega_0^\ast(S_\tau)$. Thus note that $S(Z)$ is a stack of rank-$(n-1)$ orthants with unchanged stack height
$h(S(Z))=h_S(Z\cup \{y\})$.
\\We can choose the next element $y'\in~ Y - \{y\}$, consider the projection $\pi_{y'}: \partial_y S \twoheadrightarrow
\partial_{y'}\partial_y S\cup\{\emptyset\}$. Upon putting $\pi_y(\emptyset)=\emptyset$ for all $y\in Y $, we can iterate
the argument with all elements of $Y=\{y,\ldots,z\}$, noting that only the stacks in $S({\text{Lk}}(Y))$ survive all  these projections. The composition
$$
\pi_y = \pi_z \ldots \pi_y: S \twoheadrightarrow \partial_z \ldots\partial_y \cup \{\emptyset\}
$$
projects the stacks of $S({\text{Lk}}(Y))$ onto stacks of rank-1 orthants with the original stack heights. This shows that pet$(S)$
admits a retract isomorphic to the pet-group pet$(S')$ of a disjoint union of $h(S({\text{Lk}}(Y))$ rank-1 orthants. But pet$(S')$
contains Houghton's group on $h\big(S\big({\text{Lk}}(Y)\big)\big)$ copies of $\mathbb N$ as a subgroup of finite index.
Hence Lemma 5.3 together with Ken Brown's result \cite{br87} yields  $fl\big(\text{pet}(S)\big) \leq fl\big(\text{pet}(S')\big) 
= h\big(S({\text{Lk}}(Y)\big) - 1$, as asserted.\hfill$\square$
\\
\subsection{Application to stacks of the n-skeleton of an orthant.}~ Let $S $ be a stack of $c(S)$ copies of the rank-$n$
skeleton $K^{(n)}$ of a rank-$r$ orthant $K$. The link of each cardinality-$(n-1)$ subset $Y$ of the cardinality-$r$ set
$X$ contains exactly $(r - n + 1)$ cardinality-$n$ subsets $Y'$. And $S$ contains exactly $c(S)$ orthants parallel to
$\langle Y'\rangle$.
Hence the height of disjoint union of the stacks of $S$ over the link ${\text{Lk}}(Y)$ is  $h(S({\text{Lk}}(Y)) = c(S)(r - n + 1)$.
Combining Proposition 4.2 with Theorem 5.4 thus yields

\begin{theorem} If $S$ is the rank-n skeleton of a stack of rank-r orthants then
	\[		c(S) - 1  \le  fl  \big(\emph{pet}(S)\big)  \le  c(S)(r - n + 1) - 1.\]
\end{theorem}
\begin{corollary}
If $S$ is a stack of orthants then $fl  \big(\emph{pet}(S)\big) = c(S) - 1.$\hfill $\square$
\end{corollary}

Robert Bieri\\
Binghamton University and Goethe Universität Frankfurt\\
 e-mail: bieri@math.uni-frankfurt.de ,  rbieri@math.binghamton.edu\\

Heike Sach\\
e-mail: heike\_sach@yahoo.de

\end{document}